\documentstyle[11pt]{article} \textwidth=15cm \textheight=22cm
\def\lam{\lambda}
\def\Lam{\Lambda}

\def\ol{\overline}

\def\E{{\cal E}}

\def\X{{\cal X}}
\def\Q{{\cal Q}}
\def\Y{{\cal Y}}

\def\E{{\cal E}}

\def\O{{\cal O}}
\def\R{{\cal R}}
\def\A{{\cal A}}
\def\B{{\cal B}}

\def\Q{{\cal Q}}

\def\div{{\,\rm div\,}}
\def\Term{{\,\rm Term\,}}

\def\ball{{\,\rm B\,}}
\newfont{\Blackboard}{msbm10 scaled 1200}

\newfont{\roma}{cmr10 scaled 1200}

\def\<{{\langle}}
\def\>{{\rangle}}

\def\Ga{\Gamma}
\def\var{\varphi}
\def\si{\sigma}

\def\a{\alpha}

\def\Om{\Omega}

\def\any{\forall}

\newtheorem{thm}{{}\hskip\parindent Theorem}[section]
\newtheorem{lem}{{}\hskip\parindent Lemma}[section]
\newtheorem{pro}{{}\hskip\parindent Proposition}[section]
\newtheorem{exl}{{}\hskip\parindent Example}[section]

\def\dsum{\displaystyle\sum}

\def\dfrac{\displaystyle\frac}

\def\pl{\partial}
\def\rw{\rightarrow}

\def\e{\left}
\def\i{\right}

\def\un{\upsilon}

\def\be{\begin{equation}}
\def\ee{\end{equation}}
\def\beq{\arraycolsep=1.5pt\begin{eqnarray}}
\def\eeq{\end{eqnarray}}

\def\i{\right}

\large
\title{Boundary controllability for the quasilinear wave equation}
\date{}
\author{
Peng-Fei YAO\\[0.3cm]
Key Laboratory of Control and Systems \\
Institute of Systems Science,
Academy of Mathematics and Systems Science\\
Chinese Academy of Sciences, Beijing 100080, P.R.China\\
e-mail: pfyao@iss.ac.cn\\[0.3cm]}

\begin{document}
\maketitle \footnote{This work is  supported by the NNSF of China,
grants no. 60225003, no. 60334040, and no. 60221301 }
\begin{quote}
\begin{small}
{\bf Abstract} \,\,\,We study the boundary exact controllability for
the quasilinear wave equation in the higher-dimensional case. Our
main tool is the geometric analysis. We derive the existence of long
time solutions near an equilibrium, prove the locally exact
controllability around the equilibrium under some checkable
geometrical conditions. We then establish the globally exact
controllability in such a way that the state of the quasilinear wave
equation moves from an equilibrium in one location to an equilibrium
in another location under some geometrical condition. The Dirichlet
action and the Neumann action are studied, respectively. Our results
show that exact controllability is geometrical characters of a
Riemannian metric, given by the coefficients and equilibria of the
quasilinear wave equation. A criterion of exact controllability is
given, which based on the sectional curvature of the Riemann metric.
Some examples are presented to verify the global exact
controllability.
\\[3mm]
{\bf AMS(MOS) subject classifications}\,\,\, 49B,
49E, 35B35, 35L65, 35L70, 38J45 \\[3mm]
\end{small}
\end{quote}

\def\theequation{1.\arabic{equation}}
\setcounter{equation}{0}
\section{Introduction and the main results }
\hskip\parindent

Let $\Om \subset\R^n$ be an open, bounded set with the smooth
boundary $\Ga$. Suppose that $\Ga$ consists of two disjoint parts,
$\Ga_0$ and $\Ga_1$. Let $T>0$ be given.
 We consider a controllability
problem \be\label{2.1}\cases{\ddot{u}=\sum_{ij=1}^na_{ij}(x,\nabla
u)u_{x_ix_j}+b(x,\nabla u)\quad\mbox{on}\quad (0,T)\times\Om,\cr
u=0\quad \mbox{on}\quad (0,T)\times\Ga_1, \cr u=\var\quad
\mbox{on}\quad (0,T)\times\Ga_0,\cr
u(0)=u_0,\quad\dot{u}(0)=u_1,}\ee where $a_{ij}(x,y)$, $b(x,y)$
are smooth functions on $\ol{\Om}\times\R^n$ such that
\be\label{2.2}A(x,y)=\e(a_{ij}(x,y)\i)>0\quad\any
\,\,(x,y)\in\ol{\Om}\times\R^n,\ee
\be\label{2.3}b(x,0)=0\quad\any\,\,x\in\ol{\Om}.\ee

Let $u_0$, $u_1$, $\hat{u}_0$, and $\hat{u}_1$ be given functions
on $\ol{\Om}$ and $T>0$ be given. If there is a boundary function
$\var$ on $(0,T)\times\Ga_0$ such that the solution of the problem
(\ref{2.1}) satisfies $$ u(T)=\hat{u}_0,\quad
\dot{u}(T)=\hat{u}_1\quad\mbox{on}\quad\Om,$$ we say the system
(\ref{2.1}) is exactly controllable from $(u_0,u_1)$ to
$(\hat{u}_0,\hat{u}_1)$ at time $T$ by boundary with the Dirichlet
action.

In the case of one dimension, these problems have been studied by
Cirina \cite{C}, Li and Rao \cite{LR}, Schmidt \cite{S}, and so on.
In the case of multi-dimension, $n\geq2$, very little is known in
the content of control. The work here represents a substantial
advance on this topic. The key issue is to establish the geometrical
structure of the problem: The locally exact controllability is
equivalent to the smooth control problem of a linear, variable
coefficient wave equation which is related to the geometric theory.
The detail study of the smooth control of the linear problem
provides a smooth control to the quasilinear problem. Then a
compactness principle gives the globally exact controllability. This
idea is also used to study the existence of global solutions of the
quasilinear wave equation with boundary dissipation by Yao
\cite{Y2}.

Let us choose some Sobolev spaces to formulate our problems. Let
$$m\geq[n/2]+3$$ be a given positive integer.
Inspired by Dafermos and Hrusa \cite{DH}, we assume initial data
$(u_0,u_1)\in H^m(\Om)\times H^{m-1}(\Om)$ to study the
possibility of moving it to another state in $ H^m(\Om)\times
H^{m-1}(\Om)$ at time $T$ via a boundary control $\var\in
\cap_{k=0}^{m-2}C^k\e([0,T],H^{m-k-1/2}(\Ga_0)\i)$.

In general, solutions of the system (\ref{2.1}) may below up in a
finite time even if the initial data and the boundary control are
smooth. On the other hand, in order to move one state to another,
the control time must be larger than the wave length of the
system. To cope with those situations, we shall study the locally
exact controllability of the system around an equilibrium and the
globally exact controllability form one equilibrium to another.

We say $w\in H^m(\Om)$ is an equilibrium of the system (\ref{2.1})
if \be\label{2.1*}\sum_{ij=1}^na_{ij}(x,\nabla
w)w_{x_ix_j}+b(x,\nabla w)=0\quad\mbox{on}\quad\Om.\ee

We say that $(u_0,u_1)\in H^m(\Om)\times H^{m-1}(\Om)$ and
$\var\in \cap_{k=0}^{m-2}C^k\e([0,T],H^{m-k-1/2}(\Ga_0)\i)$
satisfy the compatibility conditions of $m$ order if
\be\label{2.1****}u_k\in H^{m-k}(\Om),\quad u_k|_{\Ga_1}=0,\ee
$$\var^{(k)}(0)=u_k|_{\Ga_0},\quad k=0,\,1,\,\cdots,\,m-1,$$
where for $k\geq 2$, \be\label{2.1***}u_k=u^{(k)}(0),\ee as
computed formally (and recursively) in terms of $u_0$ and $u_1$,
using the equation in (\ref{2.1}).

Let \be\label{3.15*}H^1_{\Ga_1}(\Om)=\{\,v\,|\,v\in
H^1(\Om),\,\,v|_{\Ga_1}=0\,\}.\ee  Near one equilibrium, the
system has solutions of long time. This is the following

\begin{thm}
Let $w\in H^m(\Om)\cap H^1_{\Ga_1}(\Om)$ be an equilibrium of the
problem $(\ref{2.1})$. Let $T>0$ be arbitrary given. Then there is
$\varepsilon_T>0$, which depends on the time $T$, such that, if
$(u_0,u_1)\in \e(H^m(\Om)\cap H^1_{\Ga_1}(\Om)\i)\times
\e(H^{m-1}(\Om)\cap H^1_{\Ga_1}(\Om)\i)$ satisfy $$
\|u_0-w\|_m<\varepsilon_T,\quad \|u_1\|_{m-1}<\varepsilon_T,$$
where $\|\cdot\|_m$ denotes the norm of $H^m(\Om)$, and $\var\in
\cap_{k=0}^{m-2}C^k\e([0,T],H^{m-k-1/2}(\Ga_0)\i)$ with
$\var^{(k)}\in H^1\e((0,T)\times\Ga_0\i)$ for $0\leq k\leq m-1$
satisfies the compatibility conditions with $(u_0,u_1)$ of $m$
order and
$$\sum_{k=0}^{m-2}\|\hat{\var}^{(k)}\|^2_{C\e([0,T],H^{m-1/2-k}(\Ga_0)\i)}+
\sum_{k=0}^{m-1}\|\hat{\var}^{(k)}\|^2_{H^1\e((0,T)\times\Ga_0\i)}<\varepsilon_T,$$
where
$$\hat{\var}(t,x)=\var-w|_{\Ga_0},$$  then the system $(\ref{2.1})$ has a solution \be\label{2.1**}
u\in\cap_{k=0}^mC^k\e([0,T],H^{m-k}(\Om)\i).\ee
\end{thm}

Let $w\in H_{\Ga_1}^m(\Om)$ be an equilibrium of the system
(\ref{2.1}). We define \be\label{2.2*} g=A^{-1}(x,\nabla w)\ee as a
Riemannian metric on $\ol{\Om}$ and consider the couple
$(\ol{\Om},g)$ as a Riemannian manifold with a boundary $\Ga$. Here
the metric $g$ depends on the functions $a_{ij}(\cdot,\cdot)$ and
also on the equilibrium $w$. We denote by $\<\cdot,\cdot\>_{g}$ the
inner product induced by $g$. Let $x^0\in\ol{\Om}$ be given. We
denote by $\rho(x)=\rho(x,x^0)$ the distance function from
$x\in\ol{\Om}$ to $x^0$ under the Riemannian metric
$g$.\\

{\bf Definition}\,\,\,An equilibrium $w\in H^m(\Om)$ is called
exactly controllable if there are $x^0\in\ol{\Om}$ and $\rho_0>0$
such that\be\label{2.3*}D^2_{g}\rho^2(X,X)\geq \rho_0|X|^2_{g}\quad
\any\,X\in\Om_x,\,\,x\in\ol{\Om},\ee where $D^2_{g}\rho^2$ denotes
the Hessian of the function $\rho^2$ under the metric $g$ which is a
bilinear form on
$\ol{\Om}$.\\

The condition (\ref{2.3*}) means that the function $\rho^2(x)$ is
strictly convex on $\ol{\Om}$ under the metric $g$. This is true if
$x$ is in a neighbourhood of $x^0$. Whether it holds on the whole
domina $\ol{\Om}$ is closely related to the sectional curvature of
the Riemannian metric $g$, see some examples later. Yao \cite{Y}
presents a counterexample where the condition (\ref{2.3*}) is not
always true for all $x\in\ol{\Om}$ even when $A(x,y)=A(x)$ (the
linear problem). If the matrices $A(x,y)=A(y)$ and the equilibrium
is zero, the condition (\ref{2.3*}) holds for any $\Om\subset\R^n$
with $\rho_0=2$. A proposition below is useful to verify the
condition (\ref{2.3*}).

For $x\in\ol{\Om}$, let $\Pi\subset R^n_x$ be a two-dimensional
subspace. Denote by $k_x(\Pi)$ the sectional curvature of the
subspace $\Pi$ at $x$ under the Riemannian metric $g$.
Let\be\label{2.8*} \kappa=\sup_{x\in\ol{\Om},\,\Pi\subset
\R^n_x}k_x(\Pi).\ee

Then

\begin{pro}
If an equilibrium $w$ is such that $\kappa\leq 0$, then $w$ is
exactly controllable.

Suppose $\kappa>0$. Set
$$\lam=\inf_{x\in\ol{\Om},\,y\in\R^n,\,|y|=1}\sqrt{\<A(x,\nabla w)y,\,\,y\>}.$$
 If there is a point $x_0\in\ol{\Om}$ such that \be\label{2.34*}\ol{\Om}\subset
 B\e(x_0,\frac{\lam\pi}{2\sqrt{\kappa}}\i),\ee where $B\e(x_0,\frac{\lam\pi}{2\sqrt{\kappa}}\i)
 =\{\,x\,|\,x\in\R^n,\,|x-x_0|<\frac{\lam\pi}{2\sqrt{\kappa}}\,\}$,
then $w$ is exactly controllable.
\end{pro}

Near one equilibrium being exactly controllable, we have the
following exact controllability results:

\begin{thm}
Let an equilibrium $w\in H^m(\Om)\cap H^1_{\Ga_1}(\Om)$ be exactly
controllable. Let\be\label{2.4*}
T_0=\frac{4}{\rho_0}\sup_{x\in\ol{\Om}}\rho.\ee Furthermore, if
$\Ga_1\not=\emptyset$, we assume that\be\label{2.9*} \rho_\nu\leq
0\quad\any\,x\in\Ga_1,\ee where $\rho_\nu$ is the normal
derivative of the distance function $\rho$ of the metric $g$ with
respect to the normal $\nu$ of the dot metric of $\R^n$. Then, for
$T>T_0$ given, there is $\varepsilon_T>0$ such that, for any
$(u^i_0,u^i_1)\in H^m(\Om)\times H^{m-1}(\Om)$ with
$$
\|u^i_0-w\|_m<\varepsilon_T,\quad\|u^i_1\|_{m-1}<\varepsilon_T,\quad
u_0^i|_{\Ga_1}=u_1^i|_{\Ga_1}=0,\quad i=1,\,2,$$ we can find
$\var\in \cap_{k=0}^{m-2}C^k\e([0,T],H^{m-k-1/2}(\Ga_0)\i)$ with
$\var^{(k)}\in H^1\e((0,T)\times\Ga_0\i)$ for $0\leq k\leq m-1$
which is compatible with $(u^1_0,u^1_1)$ of $m$ order such that
the solution of the system $(\ref{2.1})$ with the initial data
$(u_0^1,u_1^1)$ satisfies\be\label{2.6*}
u(T)=u^2_0,\quad\dot{u}(T)=u^2_1.\ee
\end{thm}

The above is a local result. However, if we have enough equilibria
exactly controllable, we can move the quasilinear wave state along
a curve of equilibria, moving in successive small steps from one
equilibrium to another nearby equilibrium until the target
equilibrium is reached. This uses the open mapping theorem,
locally exact controllability,  and a compactness argument. This
approach was used by Schmidt \cite{S} for the quasilinear string.

Let $w\in H_{\Ga_1}^m(\Om)$ be a given equilibrium. For
$\a\in[0,1]$, we assume that $w_\a\in H^m(\Om)$ are the solutions
of the Dirichlet problem \be\label{2.30*}
\cases{\sum_{ij=1}^na_{ij}(x,\nabla w_\a)w_{\a x_ix_j}+b(x,\nabla
w_\a)=0\quad x\in\Om,\cr w_\a|_\Ga=\a w|_{\Ga},}\ee such that\be\label{2.31*}
\sup_{\a\in[0,1]}\|w_\a\|_m<\infty.\ee For the existence of the
classical solution to the Dirchlet problem (\ref{2.30*}), for
example, see Gilbarg and Trudinger \cite{GT}.

\begin{thm}
Let an equilibrium $w\in H^m(\Om)\cap H^1_{\Ga_1}(\Om)$ be exactly
controllable. Let $w_\a\in H^m(\Om)$, given by $(\ref{2.30*})$, be
also exactly controllable for all $\a\in[0,1]$ such that
$(\ref{2.31*})$ hold. Then, there are $T>0$ and
$$\var\in \cap_{k=0}^{m-2}C^k\e([0,T],H^{m-k-1/2}(\Ga_0)\i)$$ with
$\var^{(k)}\in H^1\e((0,T)\times\Ga_0\i)$ for $0\leq k\leq m-1$
which is compatible with the initial data $(w,0)$ such that the
solution of the system $(\ref{2.1})$ with $(u_0,u_1)=(w,0)$
satisfies
$$u(T)=\dot{u}(T)=0.$$
\end{thm}

Since the quasilinear wave equation is time-reversible, an
equilibrium can be moved to another if they can both be moved to
zero. However, this result only gives the existence of the control
time $T$. We do not know how large the $T$ is because it is
given by the compactness principle.\\

Next, we turn to the boundary control with the Neumann action. Let
$\Ga=\Ga_0\cup\Ga_1$ and $\ol{\Ga}_0\cap\ol{\Ga}_1=\emptyset$ with
$\Ga_1$ nonempty.  This time we assume that the quasilinear part
of the system is in the divergence form. Let $T>0$ be given. We
consider a controllability problem
\be\label{2.1n}\cases{\ddot{u}=\div \textbf{a}(x,\nabla
u)\quad\mbox{on}\quad (0,T)\times\Om,\cr u=0\quad \mbox{on}\quad
(0,T)\times\Ga_1, \cr \<\textbf{a}(x,\nabla u),\,\nu\>=\var\quad
\mbox{on}\quad (0,T)\times\Ga_0,\cr
u(0)=u_0,\quad\dot{u}(0)=u_1,}\ee where
$\textbf{a}(\cdot,\cdot)=\e(a_1(\cdot,\cdot),\cdots,a_n(\cdot,\cdot)\i)$
and $a_i(\cdot,\cdot)$ are smooth functions on
$\ol{\Om}\times\R^n$ such that
\be\label{2.2n}\textbf{a}(x,0)=0\quad\any\,x\in\ol{\Om};\quad
A(x,y)=\e(a_{iy_j}(x,y)\i)>0\quad\any
\,\,(x,y)\in\ol{\Om}\times\R^n.\ee In the problem (\ref{2.1n}),
$\nu$ is the normal of the boundary $\Ga$ in the dot metric of
$\R^n$.

We say $w\in H^m(\Om)\cap H^1_{\Ga_1}(\Om)$ is an equilibrium of the
system (\ref{2.1n}) if \be\label{2.1*n}\div \textbf{a}(x,\nabla
w)=0\quad\mbox{on}\quad\Om.\ee

We say that $(u_0,u_1)\in H^m(\Om)\times H^{m-1}(\Om)$ and
$\var\in \cap_{k=0}^{m-2}C^k\e([0,T],H^{m-k-3/2}(\Ga_0)\i)$
satisfy the compatibility conditions of $m$ order with the Neumann
boundary data on $\Ga_0$ and the Dirichlet data on $\Ga_1$ if
(\ref{2.1****}) hold and
\be\label{2.1**n}\var^{(k)}(0)=\cases{\<\textbf{a}(x,\nabla
u_0),\,\nu\>\quad x\in\Ga_0,\quad k=0,\cr \<A(x,\nabla u_0)\nabla
u_k,\,\nu\>\quad x\in\Ga_0,\quad 1\leq k\leq m-1,}\ee where for
$k\geq 2$, $u_k$ are given by (\ref{2.1***}).

\begin{thm}
Let $w\in H^m(\Om)\cap H^1_{\Ga_1}(\Om)$ be an equilibrium of the
problem $(\ref{2.1n})$. Let $T>0$ be arbitrary given. Then there
is $\varepsilon_T>0$, which depends on the time $T$, such that, if
$(u_0,u_1)\in \e(H^m(\Om)\cap H^1_{\Ga_1}(\Om)\i)\times
\e(H^{m-1}(\Om)\cap H^1_{\Ga_1}(\Om)\i)$ satisfy $$
\|u_0-w\|_m<\varepsilon_T,\quad \|u_1\|_{m-1}<\varepsilon_T,$$ and
$\var\in \cap_{k=0}^{m-2}C^k\e([0,T],H^{m-k-3/2}(\Ga_0)\i)$ is
such that $\var^{(k)}\in L^2\e((0,T),H^{1/2}(\Ga_0)\i)$ for $0\leq
k\leq m-1$, which satisfies the compatibility conditions
$(\ref{2.1**n})$ with $(u_0,u_1)$ of $m$ order and
$$\sum_{k=0}^{m-2}\|\hat{\var}^{(k)}\|^2_{C\e([0,T],H^{m-k-3/2}(\Ga_0)\i)}+\sum_{k=0}^{m-1}\|\hat{\var}^{(k)}
\|^2_{L^2\e((0,T),H^{1/2}(\Ga_0)\i)}<\varepsilon_T,$$ where
$$\hat{\var}(t,x)=\var-\<\textbf{a}(x,\nabla w),\,\nu\>\quad x\in\Ga_0,$$
then the system $(\ref{2.1n})$ has a solution \be\label{2.1**p}
u\in\cap_{k=0}^mC^k\e([0,T],H^{m-k}(\Om)\i).\ee
\end{thm}

\begin{thm}
Let an equilibrium $w\in H^{m+1}(\Om)\cap H^1_{\Ga_1}(\Om)$ be
exactly controllable. Let\be\label{2.9**} \rho_\nu\leq
0\quad\any\,x\in\Ga_1.\ee Then there exists a $T_0>0$ such that
the following things are true.  For any $T>T_0$ given, there is
$\varepsilon_T>0$ such that, for any $(u^i_0,u^i_1)\in
H^{m+1}(\Om)\times H^{m}(\Om)$ with
$$
\|u^i_0-w\|_{m+1}<\varepsilon_T,\quad\|u^i_1\|_{m}<\varepsilon_T,\quad
u_0^i|_{\Ga_1}=u_1^i|_{\Ga_1}=0,\quad i=1,\,2,$$ we can find
$\var\in \cap_{k=0}^{m-2}C^k\e([0,T],H^{m-k-3/2}(\Ga_0)\i)$ with
$\var^{(k)}\in L^2\e((0,T),H^{1/2}(\Ga_0)\i)$ for $0\leq k\leq
m-1$ which is compatible with $(u^1_0,u^1_1)$ of $m$ order such
that the solution of the system $(\ref{2.1n})$ with the initial
data $(u_0^1,u_1^1)$ satisfies\be\label{2.6**}
u(T)=u^2_0,\quad\dot{u}(T)=u^2_1.\ee
\end{thm}

Here we lose an explicit formula of $T_0$.

Unlike the control with the Dirichlet action, we only have the exact
controllability results in the space $H^{m+1}(\Om)\times H^{m}(\Om)$
by a control $\var\in
\cap_{k=0}^{m-2}C^k\e([0,T],H^{m-k-3/2}(\Ga_0)\i)$ with
$\var^{(k)}\in L^2\e((0,T),H^{1/2}(\Ga_0)\i)$ for $0\leq k\leq m-1$.
This is because the Neumann action loses a regularity of $1$ order
(actually, $1/2$ order), see Theorem 2.2 in the end of Section 2. In
addition, although we can move one state to another in the space
$H^{m+1}(\Om)\times H^{m}(\Om)$, we can not guarantee the solution
$(u(t),\dot{u}(t))$ of the problem (\ref{2.1n}) always stays in
$H^{m+1}(\Om)\times H^{m}(\Om)$ in the process of the control for
$0\leq t\leq T$ where they are actually in the space $H^m(\Om)\times
H^{m-1}(\Om)$ for all $t\in [0,T]$ by Theorem 1.4. The same things
happen to the globally exact controllability results in Theorem 1.6
below.

Let an equilibrium $w\in H_{\Ga_1}^{m+1}(\Om)$ be given. For $\a\in[0,1]$, we
assume that $w_\a\in H^{m+1}(\Om)$ are the solutions of the
Dirichlet problem (\ref{2.30*}) with, this time, an uniform bound
\be\label{n} \sup_{\a\in[0,1]}\|w_\a\|_{m+1}<\infty.\ee

\begin{thm}
Let an equilibrium $w\in H^{m+1}(\Om)\cap H^1_{\Ga_1}(\Om)$ be exactly
controllable. Let $w_\a\in H^{m+1}(\Om)$ be also exactly
controllable for all $\a\in[0,1]$ such that $(\ref{n})$ hold. Then,
there are $T>0$ and
$$\var\in \cap_{k=0}^{m-2}C^k\e([0,T],H^{m-k-3/2}(\Ga_0)\i)$$
 with
$\var^{(k)}\in L^2\e((0,T),H^{1/2}(\Ga_0)\i)$ for $0\leq k\leq m-1$,
which is compatible with the initial data $(w,0)$ of $m$ order such
that the solution of the system $(\ref{2.1n})$ with
$(u_0,u_1)=(w,0)$ satisfies
$$u(T)=\dot{u}(T)=0.$$
\end{thm}

Boundary exact controllability on linear problems has been
developing since 70's and very active in recent years. We mention
Bardos, Lebeau, Rauch \cite{BLR}, Castro, Zuazua \cite{CZ}, Egorov
\cite{E}, Fattorini \cite{Fa}, Ho \cite{H}, Lasiecka, Triggiani
\cite{LT}, Lions \cite{L}, Russel \cite{Ru},
Seidman \cite{Se}, Tataru \cite{Ta},  Yao \cite{Y}, \cite{Y1}, Yong, Zhang \cite{YZ}, just a few.\\

Finally, let us see some examples to verify Theorem 1.3.

\begin{exl}
Let $n=2$ and $m=4$. Consider the control problem\be\label{2.10*}
\cases{\ddot{u}=(|\nabla u|^2+1)\Delta u\quad
(t,x)\in(0,T)\times\Ga,\cr u|_{\Ga}=\var\quad 0\leq t\leq T, \cr
u(0)=w\quad\dot{u}(0)=0\quad  x\in\Om,}\ee where
$\Delta=\dfrac{\pl^2}{\pl x^2}+\dfrac{\pl^2}{\pl y^2}$.

Then $w\in H^4(\Om)$ is an equilibrium if and only if
$$\Delta w=0\quad x\in\Om.$$ Let $w\in H^4(\Om)$ be an equilibrium.
Then metric $(\ref{2.2*})$ is given by
$$g=A^{-1}(x,\nabla w),\quad A(x,\nabla w)=\e(\begin{array}{cc}|\nabla w|^2+1&0\\
0&|\nabla w|^2+1\end{array}\i).$$

By Lemma 3.2, Yao $\cite{Y}$, the Gauss curvature of the Riemmannin
manifold $(\ol{\Om},g)$ is
$$k(x)=|D^2w|^2,\quad x\in\ol{\Om},$$ and
$$\kappa=\sup_{x\in\ol{\Om}}|D^2w|^2,$$ where $D^2w$ is the Hessian of $w$ in the dot metric of
$\R^2$. {\bf Then the zero equilibrium, $w=0$, is exactly
controllable for any $\Om\subset\R^n$}. In addition, we have the
conclusion: {\bf If two equilibria $w_i\not=0$ in $H^4(\Om)$ are
such that there are $x_i\in\ol{\Om}$ satisfying
\be\label{2.13**}\ol{\Om}\subset B(x_i,\gamma_i),\ee where
\be\label{2.13*}\gamma_i=\frac{\pi\sqrt{1+\inf_{x\in\ol{\Om}}|\nabla
w_i|^2}}{2\sup_{x\in\ol{\Om}}|D^2w_i|},\ee for $i=1$, $2$, then
there are a control time $T>0$ and a control function $$\var\in
\cap_{k=0}^2C^k\e([0,T],H^{7/2-k}(\Ga)\i)$$ with $\var^{(k)}\in
H^1\e((0,T)\times\Ga\i)$ for $0\leq k\leq 3$ such that the solution
of the problem $(\ref{2.10*})$ with the initial $(w_1,0)$ satisfies
$$u(T)=w_2,\quad\dot{u}(T)=0.$$}
Suppose that the two equilibria $w_i$ in $H^4(\Om)$ are such that
the conditions $(\ref{2.13**})$ are true.
 Then, for $\a\in[0,1]$, $w_{i\a}=\a w_i$ are equilibria with
$w_{i\a}|_\Ga=\a w_i|_\Ga$. Since
$$\gamma_i\leq
\frac{\pi\sqrt{1+\a^2\inf_{x\in\ol{\Om}}|\nabla
w_i|^2}}{2\a\sup_{x\in\ol{\Om}}|D^2w_i|},$$ for all $0<\a\leq1$, the
conditions $(\ref{2.13**})$ are true for all $w_{i\a}$ with
$a\in[0,1]$. By Theorem 1.3 and Proposition 1.1, the initial data
$(w_i,0)$ can be moved to $(0,0)$, respectively.

Let $$w_1=a(x^2-y^2),\quad w_2=axy,\quad \Om=\mbox{the unit disc},$$
where $0<a<1/(2\sqrt{2})$. It is easy to check that $w_i$ meet the
conditions $(\ref{2.13**})$ for $i=1$, $2$. Then the state of the
system $(\ref{2.10*})$ can be moved from $(w_1,0)$ to $(w_2,0)$ at
some time $T>0$.
\end{exl}

\begin{exl}
Consider the control problem\be\label{2.7*}
\cases{\ddot{u}=(|\nabla u|^2+1)^{-1}\Delta u\quad
(t,x)\in(0,T)\times\Ga,\cr u|_{\Ga}=\var\quad 0\leq t\leq T, \cr
u(0)=w\quad\dot{u}(0)=0\quad x\in\Om.}\ee

Let $w\in H^4(\Om)$ be an equilibrium. The metric is
$$g=\e(\begin{array}{cc}|\nabla w|^2+1&0\\
0&|\nabla w|^2+1\end{array}\i).$$ The Gauss curvature of $(\ol{\Om},
g)$ is $$ k(x)=-\frac{|D^2w|^2}{(|\nabla w|^2+1)^3}\leq
0,\quad\any\,\,x\in\ol{\Om}.
$$ We have the conclusion: {\bf For any two equilibria $w_1$,
$w_2\in H^4(\Om)$ and any $\Om\subset\R^n$, there are a control time
$T>0$ and a control function $\var\in
\cap_{k=0}^2C^k\e([0,T],H^{7/2-k}(\Ga)\i)$ with $\var^{(k)}\in
H^1\e((0,T)\times\Ga\i)$ for $0\leq k\leq 3$ such that the state of
the system $(\ref{2.7*})$ is moved from $(w_1,0)$ to $(w_2,0)$}.
\end{exl}

\def\theequation{2.\arabic{equation}}
\setcounter{equation}{0}
\section{Solutions of long time }
\hskip\parindent The basic results of the existence of short time
solutions to the quasilinear wave equation has been established by
Dafermos and Hrusa \cite{DH}. We here only study some energy
estimates of the short time solutions to have long time solutions
when initial data are close to an equilibrium.

Let $(u_0,u_1)\in H^m(\Om)\times H^{m-1}(\Om)$ and $\var\in
C^{\infty}\e((0,T)\times\Ga_0\i)$ satisfy the compatibility
conditions of $m$ order. If we extend $\var$ from $(0,T)\times\Ga_0$
to $(0,T)\times\Om$, still denoted by $\var$ and let
\be\label{1.32}v=u-\var\ee as a new unknown, then the problem will
have solutions $v\in\cap_{k=0}^mC^k\e([0,T],H^{m-k}(\Om)\i)$ of
short time by Dafermos and Hrusa \cite{DH}, Theorem 5.1.

To obtain solutions of long time near an equilibrium location, we
need to estimate the energy of solutions to the problem (\ref{2.1}).
We observe that, if we apply Dafermos and Hrusa \cite{DH}, Theorem
3.1 to our problem  after the transform (\ref{1.32}), we shall see
that the regularity of
$\var\in\cap_{k=0}^{m-1}C^k\e([0,T],H^{m-k-1/2}(\Ga_0)\i)$ is
insufficient to guarantee
$u\in\cap_{k=0}^mC^k\e([0,T],H^{m-k}(\Om)\i)$ because we have lost a
regularity of $1/2$ order by the transform (\ref{1.32}). For this
reason, we shall here work out our energy estimates starting from
the problem (\ref{2.1}) directly.

We suppose that the equilibrium is the zero, $w=0$, in this section.
If an equilibrium $w\in H^m(\Om)\cap H^1_{\Ga_1}(\Om)$ is not zero,
we can make an transform by
$$u=w+v,$$ and consider the $v$-problem
$$\cases{\ddot{v}=\sum_{ij}\hat{a}_{ij}(x,\nabla
v)v_{x_ix_j}+\hat{b}(x,\nabla v)\quad (t,x)\in(0,T)\times\Om,\cr
v|_{\Ga_1}=0,\quad v|_{\Ga_0}=\var-w|_{\Ga_0},\cr v=v_0,\quad
\dot{v}(0)=v_1,}$$ where $$\hat{a}_{ij}(x,y)=a_{ij}(x,\nabla w+y),$$
$$\hat{b}(x,y)=\sum_{ij}a_{ij}(x,\nabla w+y)w_{x_ix_j}+b(x,\nabla
w+y),$$
$$v_0=w_0-w,\quad v_1=w_1.$$

Let $u\in\cap_{k=0}^mC^k\e([0,T],H^{m-k}(\Om)\i)$ be a solution of
the problem (\ref{2.1}) for some $T>0$. Suppose that
$$\var\in\cap_{k=0}^{m-2}C^k\e([0,T],H^{m-k-1/2}(\Ga_0)\i)$$
and $$\var^{(k)}\in H^1\e((0,T)\times\Ga_0\i),\quad
0\leq k\leq m-1.$$ We introduce
\be\label{1.28}\E(t)=\sum_{k=0}^m\|u^{(k)}(t)\|^2_{m-k},\quad
\E_\Ga(t)=\sum_{k=0}^{m-2}\|\var^{(k)}(t)\|^2_{m-k-1/2,\Ga_0},\ee
\be\label{1.30}Q(t)=\sum_{k=1}^m\e(\|u^{(k)}(t)\|^2+\|\nabla
u^{(k-1)}(t)\|^2\i),\ee \be\label{1.31}
Q_\Ga(t)=\sum_{k=1}^{m}\e(\|\var^{(k)}(t)\|_{\Ga_0}^2+\|\nabla
\var^{(k-1)}(t)\|_{\Ga_0}^2\i),\ee where $\|\cdot\|_0=\|\cdot\|$ and
$\|\cdot\|_\Ga=\|\cdot\|_{0,\Ga_0}$ are norms of $L^2(\Om)$ and
$L^2(\Ga_0)$ and $\|\cdot\|_j$, $\|\cdot\|_{j,\Ga_0}$ are norms of
$H^j(\Om)$, $H^j(\Ga_0)$ for $1\leq j\leq m$, respectively.

\begin{thm}\label{t2.1}
We consider solutions of the problem $(\ref{2.1})$ near the zero
equilibrium. Let $\gamma>0$ be given and $u$ be a solution of the
problem $(\ref{2.1})$ on the interval $[0,T]$ for some $T>0$ such
that \be\label{1.25} \sup_{0\leq t\leq T}\|u(t)\|_m\leq \gamma.\ee
Then there is $c_\gamma>0$, which depends on the $\gamma$ and but
is independent of initial data $(u_0,u_1)$ and boundary functions
$\var$, such that \be\label{2.26} Q(t)\leq\E(t)\leq c_\gamma
Q(t)+c_\gamma\E_\Ga(t)+c_\gamma\sum_{k=2}^m\E^{k}(t),\quad 0\leq
t\leq T,\ee and \be\label{1.12} Q(t)\leq c_\gamma Q(0)
+c_\gamma\int_0^t\e[\e(1+\E^{1/2}(t)\i)Q(t)+Q_\Ga(t)+\sum_{k=2}^m\E^k(t)\i]dt,\ee
for $t\in[0,T]$.
\end{thm}

We collect here a few basic properties of Sobolev spaces to be
invoked in the sequel.

(i) Let $s_1>s_2\geq0$. For any $\varepsilon>0$ there is
$c_\varepsilon>0$ such that\be\label{2.29} \|w\|^2_{s_2}\leq
\varepsilon\|w\|^2_{s_1}+c_\varepsilon\|w\|^2\quad\any\,\,w\in
H^{s_1}(\Om).\ee

(ii) If $s>n/2$, then for each $k=0$, $\cdots$, we have
$H^{s+k}(\Om)\subset C^k(\ol{\Om})$ with continuous inclusion.

(iii) If $r$:$=\min\{s_1,s_2,s_1+s_2-[n/2]-1\}\geq0$, then there is
a constant $c>0$ such that\be\label{2.6} \|fg\|_r\leq
c\|f\|_{s_1}\|g\|_{s_2}\quad\any\,\,f\in H^{s_1}(\Om),\,\,g\in
H^{s_2}(\Om).\ee

(iv) Let $s_j\geq0$, $j=1$, $\cdots$, $k$, and $r$:$=\min_{1\leq
i\leq k}\min_{j_1\leq\cdots\leq
j_i}\{s_{j_1}+\cdots+s_{j_i}-(i-1)([n/2]+1)\}\geq0$. Then there is
a constant $c>0$ such that\be\label{2.7} \|f_1\cdots f_k\|_r\leq
c\|f_1\|_{s_1}\cdots\|f_k\|_{s_k}\quad\any\,\,f_j\in
H^{s_j}(\Om),\,\,1\leq j\leq k.\ee

Let $u\in\cap_{k=0}^mC^k\e([0,T],H^{m-k}(\Om)\i)$ be a solution of
short time to the problem (\ref{2.1}). We introduce a linear
operator $B(t)$ by $$B(t)w=-\sum_{ij=1}^na_{ij}(x,\nabla
u)w_{x_ix_j}\quad w\in H^1(\Om).$$ Then\be\label{1.3}
(B(t)w,v)=-\int_\Ga v\<A\nabla w,\nu\>d\Ga+(A\nabla w,\nabla
v)-(Cw,v),\quad w,\,\,v\in H^1(\Om),\ee where
$A=\e(a_{ij}(x,\nabla u)\i)$, $\nu$ is the normal of $\Ga$ in the
dot metric, and
$$Cw=-\sum_{ij=1}^n\e(a_{ij}(x,\nabla u)\i)_{x_j}w_{x_i}.$$
Then the problem (\ref{2.1}) becomes\be\label{1.13}
\cases{\ddot{u}(t)+B(t)u(t)=b(x,\nabla u)\quad
(t,x)\in(0,T)\times\Om, \cr u(t)|_{\Ga_1}=0,\quad
u(t)|_{\Ga_0}=\var(t),\quad 0\in (0,T),\cr u(0)=w_0,\quad
\dot{u}(0)=w_1,\quad x\in\Om.}\ee

\begin{lem}

$(i)$ Let $f(x,y)$ be a smooth function on $\ol{\Om}\times\R^n$.
Set $F(x)=f(x,\nabla u)$. For $0\leq k\leq m-1$, there is
$c=c(\sup_{x\in\Om}|\nabla u|)>0$ such that\be\label{2.24}
\|F\|_k\leq c\sum_{j=0}^k(1+\|u\|_m)^j.\ee

$(ii)$ Let $u$ be a solution of the problem $(\ref{2.1})$ and
$\gamma>0$ be given. Suppose that the condition $(\ref{1.25})$
holds true. Then there is $c_\gamma>0$, which depends on the
$\gamma$, such that\be\label{1.22} \|v\|^2_{k+1}\leq
c_\gamma\e(\|B(t)v\|^2_{k-1}+\|v\|^2_{k+1/2,\,\Ga}+\|v\|^2_{k}\i),\quad
v\in H^k(\Om),\ee for $0\leq k\leq m-1$.

\end{lem}

{\bf Proof.}\,\,\,(i) By induction. The inequality (\ref{2.24}) is clearly true for
$k=0$. Suppose that it holds for $0\leq k<m-1$. Since
$$F_{x_i}=f_{x_i}(x,\nabla u)+\sum_{j=1}^nf_{y_j}(x,\nabla
u)u_{x_ix_j},\quad 1\leq i\leq n,$$ by using the formula
(\ref{2.6}) and the induction assumption for $f_{x_i}(x,\nabla
u)$ and for $f_{y_j}(x,\nabla u)$, respectively, we obtain\beq
\|F\|_{k+1}&&=(\|F\|^2+\sum_{i=1}^n\|F_{x_i}\|^2_{k})^{1/2}\nonumber\\
&&\leq c+c\sum_{i=1}^n\|f_{x_i}(x,\nabla
u)\|_{k}+c\sum_{ij=1}^n\|f_{y_j}(x,\nabla
u)\|_{k}\|u_{x_ix_j}\|_{m-2}\nonumber\\
&&\leq
c+c\sum_{j=0}^k(1+\|u\|_{m})^j+c\sum_{j=0}^k(1+\|u\|_{m})^j\|u\|_{m}\nonumber\\
&&\leq c\sum_{j=0}^{k+1}(1+\|u\|_{m})^j.\nonumber \eeq

(ii) A standard method as to the linearly elliptic problem can
give the inequality (\ref{1.22}), for example see Taylor \cite{T}.

\begin{lem}
Let $\gamma>0$ be given and $u$ be a solution of the problem
$(\ref{2.1})$ on the interval $[0,T]$ for some $T>0$ such that the
condition $(\ref{1.25})$ holds true. Then there is $c_\gamma>0$,
which depends on the $\gamma$, such that\be\label{1.16}
\|b^{(k)}(x,\nabla u)\|^2_{m-k-2}\leq c_\gamma\|
u^{(k)}(t)\|^2_{m-k-1}+ c_\gamma\sum_{i=2}^{k}\E^{i}(t),\quad
1\leq k\leq m-1,\ee and \be\label{1.17}
\|B^{(j)}(t)u^{(k-j)}(t)\|^2_{m-k-2}\leq
c_\gamma\sum_{i=1}^k\E^{1+i}(t),\quad 1\leq j\leq k\leq m-1.\ee
\end{lem}

{\bf Proof.}\,\,\,We have\be\label{1.18} b^{(k)}(x,\nabla
u)=\sum_{i=1}^k\dsum_{r_1+\cdots+r_i=k}D_y^ib\e(\nabla
u^{(r_1)}(t),\cdots,\nabla u^{(r_i)}(t)\i),\ee where $D_y^ib$
denotes the covariant differential of $i$ order of the function
$b(x,y)$ with respect to the variable $y$ in the dot metric of
$\R^n$.

Let us see the term $D_yb(\nabla u^{(k)})$ first. We have $$
 D_yb(\nabla u^{(k)})=\sum_{l=1}^nb_{y_l}(x,\nabla
 u)u^{(k)}_{x_l}(t).$$ By (\ref{2.6}) and (\ref{2.24})
 \be\label{1.26}\|b_{y_l}u_{x_l}^{(k)}\|_{m-k-2}\leq c
 \|b_{y_l}\|_{m-1}\|u^{(k)}(t)\|_{m-k-1}\leq
 c_\gamma\|u^{(k)}(t)\|_{m-k-1}.\ee

Let $2\leq i\leq k$. We observe that $D_y^ib\e(\nabla
u^{(r_1)}(t),\cdots,\nabla u^{(r_i)}(t)\i)$ is a sum of terms such
as $$ f(x,\nabla u)u_{x_{l_1}}^{(r_1)}(t)\cdots
u_{x_{l_i}}^{(r_i)}(t)$$ where $r_1+\cdots+r_i=k$. Using (\ref{2.7}) and (\ref{2.24}), we
have \be\label{1.27}\|fu_{x_{l_1}}^{(r_1)}\cdots
u_{x_{l_i}}^{(r_i)}\|_{m-k-2}\leq
c_\gamma\|u^{(r_1)}(t)\|_{m-r_1}\cdots\|u^{(r_i)}(t)\|_{m-r_i}\leq
c_\gamma\E^{i/2}(t).\ee

(\ref{1.16}) follows from (\ref{1.26}) and (\ref{1.27}).  A
similar argument yields (\ref{1.17}).

\begin{lem}
Let $\gamma>0$ be given and $u$ be a solution of the problem
$(\ref{2.1})$ on the interval $[0,T]$ for some $T>0$ such that the
condition $(\ref{1.25})$ holds true. Then there is $c_\gamma>0$ such that the inequality
$(\ref{2.26})$ is true.
\end{lem}

{\bf Proof.}\,\,\,It is clear that \be\label{2.27}
\|u^{(m)}(t)\|^2+\|u^{(m-1)}(t)\|^2_1=\|u^{(m)}(t)\|^2+\|\nabla
u^{(m-1)}(t)\|^2+\|u^{(m-1)}(t)\|^2\leq \Q(t).\ee Proceeding by
induction, we assume that for some $1\leq j\leq m$ \be\label{2.28}
\|u^{(j)}(t)\|^2_{m-j}\leq c_\gamma
Q(t)+c_\gamma\E_\Ga(t)+c_\gamma\sum_{k=2}^m\E^{k}(t),\ee which, as
shown above, is true for $j=m$ and $j=m-1$. Formal differentiation
of the equation in (\ref{1.13}) $j-2$ times with respect to $t$
yields \be\label{2.30}
u^{(j)}(t)+B(t)u^{(j-2)}(t)=b^{(j-2)}(x,\nabla
u)-\sum_{k=1}^{j-2}C_kB^{(k)}(t)u^{(j-2-k)}(t).\ee Using Lemmas 2.1,
2.2, (\ref{2.30}), and (\ref{2.29}), we obtain
\beq\label{2.31}\|u^{(j-2)}(t)\|^2_{m-j+2}&&\leq
c_\gamma\|B(t)u^{(j-2)}(t)\|^2_{m-j}+c_\gamma\|\var^{(j-2)}(t)\|^2_{m-j+3/2,\,\Ga_0}+c_\gamma\|u^{(j-2)}(t)\|^2_{m-j+1}\nonumber\\
&&\leq
c_\gamma\|u^{(j)}(t)\|^2_{m-j}+c_\gamma\|u^{(j-2)}(t)\|^2_{m-j+1}+c_\gamma\E_\Ga(t)+c_\gamma\sum_{i=2}^{j-1}\E^i(t).
\eeq The inequality (\ref{2.26}) follows by induction by using the
following inequality in (\ref{2.31}) $$\|u^{(j-2)}(t)\|^2_{m-k+1}
\leq
\varepsilon\|u^{(j-2)}(t)\|^2_{m-j+2}+c_{\gamma,\varepsilon}\|u^{(j-2)}(t)\|^2.$$

\begin{lem}
Let $\gamma>0$ be given and $u$ be a solution of the problem
$(\ref{2.1})$ on the interval $[0,T]$ for some $T>0$ such that the
condition $(\ref{1.25})$ holds true. Let $w\in
H^1\e((0,T)\times\Om\i)$ solve the linear
problem\be\label{1.4*}\cases{\ddot{w}(t)+B(t)w(t)=F(t)\quad (t,x)\in
(0,T)\times\Om,\cr w|_{\Ga_1}=0,\quad w|_{\Ga_0}=\var\quad
(t,x)\in(0,T)\times\Ga,\cr w(0)=w^0,\quad\dot{w}(0)=w^1,\quad
x\in\Om.}\ee Set$$ \Upsilon(t)=\|\dot{w}(t)\|^2+\|\nabla
w(t)\|^2,\quad
\Upsilon_\Ga(t)=\|\dot{\var}(t)\|^2_{\Ga_0}+\|\nabla_\Ga\var\|^2_{\Ga_0},$$
where $\nabla_\Ga$ is the gradient of $\Ga$ in the induced metric by
the dot metric of $\R^n$. Then there is $c_\gamma>0$ such
that\be\label{1.5} \Upsilon(t)\leq c_\gamma
\Upsilon(0)+c_\gamma\int_0^t\e[(1+\|\dot{u}(t)\|_{m-1})\Upsilon(\tau)+\Upsilon_\Ga(\tau)+\|F(\tau)\|^2\i]d\tau\quad
0\leq t\leq T.\ee
\end{lem}

{\bf Proof.}\,\,\,Let
\be\label{1.6n}P(t)=\|\dot{w}(t)\|^2+\e(A\nabla w,\nabla w\i).\ee
Using (\ref{1.3}) and (\ref{1.4}), we obtain\beq\label{1.6}
\dot{P}(t)&&=2\e(\ddot{w}(t),\dot{w}(t)\i)+2\e(A\nabla
w(t),\nabla\dot{w}(t)\i)+\e(\dot{A}\nabla w,\nabla w\i)\nonumber\\
&&=2\e(F+Cw,\,\,\dot{w}\i)+\e(\dot{A}\nabla w,\nabla w\i)+
2\int_{\Ga_0}\dot{\var}w_{\nu_A}d\Ga,\eeq where
$w_{\nu_A}=\<A(x,\nabla u)\nabla w,\,\,\nu\>$. It follows from
(\ref{1.6}) that\beq\label{1.7} \Upsilon(t) &&\leq c_\gamma P(t)\leq
c_\gamma P(0)+c_\gamma\int_0^t\e[(\|\nabla
w\|+\|F(t)\|)\|\dot{w}\|+\|\dot{u}\|_{m-1}\|\nabla w\|^2\i]dt\nonumber\\
&&\quad+\varepsilon\int_0^t\int_{\Ga_0} w^2_{\nu_A}d\Ga
dt+c_{\gamma,\varepsilon}\int_0^t\int_{\Ga_0}\dot{\var}^2d\Ga dt,
\quad 0\leq t\leq T,\eeq where $\varepsilon>0$ will be determined
later.

To obtain (\ref{1.5}) from (\ref{1.7}), we have to estimate the
term $\int_0^t\int_\Ga w^2_{\nu_A}d\Ga dt$.

We now introduce a Riemannian metric $$g=A^{-1}(x,\nabla u)$$ on
$\ol{\Om}$ so that the couple $(\ol{\Om},g)$ is a Riemannian
manifold. Let $H$ be a vector field on $\ol{\Om}$ such that
\be\label{1.8}H|_{\Ga_1}=0,\quad \mbox{in a neigborhood of
$\Ga_1$};\quad H|_{\Ga_0}=\nu_A.\ee We have the following formula
(see Yao \cite{Y}, Lemma 2.1)\beq\label{1.9} \<A\nabla
w,\nabla\e(H(w)\i)\>&&=D_{g}H\e(\nabla_{g}w,\nabla_{g}w\i)+\frac{1}{2}\div\e(|\nabla_{g}w|_{g}^2H\i)\nonumber\\
&&\quad-\frac{1}{2}|\nabla_{g}w|_{g}^2\div H,\eeq where $D_{g}H$ is
the covariant differential of the vector field $H$ and
$\nabla_{g}=A(x,\nabla u)\nabla$ is the gradient of the Riemmannian
metric $g$.

We multiply the two sides of the equation in (\ref{1.4*}) by $H(w)$
and integrate over $\Om$ by parts, via the formulas (\ref{1.3}),
(\ref{1.8}), and (\ref{1.9}) to obtain
\beq\label{1.10}&&\int_{\Ga_0}\e[w_{\nu_A}^2+\frac{1}{2}(\dot{\var}^2-|\nabla_{g}w|_{g}^2)|\nu_A|^2_{g}\i]d\Ga\nonumber\\
&&=\frac{\pl}{\pl
t}\e(\dot{w},H(w)\i)+\int_\Om\e[D_{g}H(\nabla_{g}w,\nabla_{g}w)+\frac{1}{2}(\dot{w}^2-|\nabla_{g}w|_{g}^2
)\div H\i]dx\nonumber\\
&&\quad-\int_\Om\e(Cw+F(t)\i)H(w)dx. \eeq Using the formula
(\ref{1.10}) and the relation
$$|\nabla_{g}w|_{g}^2=\frac{w_{\nu_A}^2}{|\nu_A|^2_{g}}+|\nabla_{\Ga_{g}}\var|_{g}^2,\quad
x\in\Ga,$$ where $\nabla_{\Ga_{g}}$ is the gradient of $\Ga$ in the
induced metric by the Riemannian metric $g$,  we
have\be\label{1.11}\int_0^t\int_{\Ga_0} w_{\nu_A}^2d\Ga dt\leq
c_\gamma\e[\Upsilon(t)+\Upsilon(0)\i]+c_\gamma\int_0^t\int_\Om\e[\Upsilon(t)+\Upsilon_\Ga(t)\i]dxdt.
\ee

Finally, we insert the inequality (\ref{1.11}) into the inequality
(\ref{1.7}), choose a $\varepsilon>0$ so small such that the term
$\varepsilon c_\gamma \Upsilon(t)$ can be moved to the left hand
side of the inequality to obtain the inequality (\ref{1.5}).
{\bf\textbf{$\|$}}\\

{\bf The Proof of Theorem 2.1}\,\,\,Lemma 2.3 gives the inequality
(\ref{2.26}). Let us prove the inequality (\ref{1.12}).

We take $w=u^{(j-2)}(t)$ for $2\leq j\leq m+1$ from the equation
(\ref{2.30}) and apply Lemma 2.4 to obtain\beq\label{1.35}
&&\|u^{(j-1)}(t)\|^2+\|\nabla u^{(j-2)}(t)\|^2\nonumber\\
&&\leq c_\gamma
Q(0)+c_\gamma\int_0^t\e[\e(1+\E^{1/2}(t)\i)Q(t)+Q_\Ga(t)\i]dt\nonumber\\
&&\quad
c_\gamma\int_0^t\e(\|b^{(j-2)}\|^2+\sum_{k=1}^{j-2}\|B^{(k)}(t)u^{(j-2-k)}(t)\|^2\i)dt.\eeq

In addition, a similar computation as in Lemma 2.2 yields
\beq\label{1.36}\|b^{(j-2)}\|^2&&\leq c_\gamma\|\nabla
u^{(j-2)}\|^2+c_\gamma\sum_{k=2}^m\E^k(t)\nonumber\\
&&\leq c_\gamma Q(t)+c_\gamma\sum_{k=2}^m\E^k(t),\eeq
and\be\label{1.37} \|B^{(k)}(t)u^{(j-2-k)}(t)\|^2\leq
c_\gamma\sum_{k=2}^m\E^k(t),\quad 1\leq k\leq j-2.\ee

The inequality (\ref{1.12}) follows from
(\ref{1.35})-(\ref{1.37}). {\bf\textbf{$\|$}}\\

{\bf The Proof of Theorem 1.1}\,\,\,Clearly, it will suffice to
prove Theorem 1.1 for the zero equilibrium $w=0$.

Let $T_1>0$ be arbitrary given. We take $\gamma=1$. Let
\be\label{1.44}c_1=c_\gamma\geq 1\ee be fixed such that the
corresponding inequalities (\ref{2.26}) and (\ref{1.12}) of
Theorem 2.1 hold for $t$ in the existence interval of the solution
$u$, respectively.

We shall prove that, if initial data $(w_0,w_1)$ and boundary value
$\var$ are compatible of $m$ order to satisfy\be\label{1.45}
\E(0)+\max_{0\leq t\leq
T_1}\E_\Ga(t)+\int_0^{T_1}Q_\Ga(t)dt\leq\frac{1}{16c_1^3}e^{-4c_1^2T_1},\ee
then the solution of the problem (\ref{2.1}) exists at least on the
interval $[0,T_1]$.

We set
\be\label{1.43}\eta=\frac{1}{4c_1}\leq\frac{1}{4}<\frac{1}{2}.\ee
Since $\E(0)\leq\eta/4$, the solution of short time must satisfy
\be\label{1.42}\E(t)\leq\eta\leq 1/2\ee for some interval
$[0,\delta]$.

Let $\delta_0$ be the largest number such that (\ref{1.42}) is
true for $t\in[0,\delta_0)$. We shall prove $\delta_0\geq T_1$ by
contradiction.

Suppose that $\delta_0<T_1$. In this interval $[0,\delta_0]$ the
condition (\ref{1.25}) is true, we apply Theorem 2.1, and the
inequalities (\ref{2.26}) and (\ref{1.12}), via (\ref{1.44}),
(\ref{1.43}), and (\ref{1.42}), imply\be\label{1.41}
\E(t)\leq2c^2_1\e[\E(0)+\max_{0\leq t\leq
T_1}\E_\Ga(t)+\int_0^{T_1}Q_\Ga(t)dt\i]+4c_1^2\int_0^t\E(t)dt,\ee
for $t\in[0,\delta_0]$. By (\ref{1.45}) and (\ref{1.41}), the
Gronwall inequality yields $$\E(\delta_0)\leq\eta/2<\eta.$$ This is
a contradiction. {\bf\textbf{$\|$}}\\

{\bf Proof Theorem 1.4}\,\,\,This proof follows by an similar
argument in the proof of Theorem 1.1 which this time is based on
the
estimates of the following theorem.\\

We turn to the problem (\ref{2.1n}) with the Neumann data on the
portion $\Ga_0$ of the boundary. Let
$u\in\cap_{k=0}^mC^k\e([0,T],H^{m-k}(\Om)\i)$ be a solution of the
problem (\ref{2.1n}) for some $T>0$. We introduce an operator
\be\label{xn1} \B(t)v=\div B\nabla v\quad v\in H^2(\Om),\ee where
$$B=\e(b_{ij}(x,\nabla u)\i),\quad
b_{ij}(x,y)=\int_0^1a_{iy_j}(x,\si y)d\si.$$ Then \be\label{xn2}
\div \textbf{a}(x,\nabla u)=\div B\nabla u,\quad
\var=\<\textbf{a}(x,\nabla u),\,\nu\>=\<B(x,\nabla u)\nabla
u,\,\nu\>=u_{\un_B}.\ee

Suppose that
\be\label{1.38n}\var\in\cap_{k=0}^{m-2}C^k\e([0,T],H^{m-k-3/2}(\Ga_0)\i).\ee
We introduce $$\E_{\Ga
N}(t)=\sum_{k=0}^{m-2}\|\var^{(k)}(t)\|^2_{m-k-3/2,\Ga_0},$$
$$Q_{\Ga
N}=\sum_{k=0}^{m-1}\|\var^{(k)}\|^2_{1/2,\Ga_0}dt.$$

Then

\begin{thm}\label{t2.2}
Let $\gamma>0$ be given and $u$ be a solution of the problem
$(\ref{2.1n})$ on the interval $[0,T]$ for some $T>0$ such that
the inequality $(\ref{1.25})$ is true.  Then there is
$c_\gamma>0$, which only depends on the $\gamma$, such that
\be\label{2.26n} Q(t)\leq\E(t)\leq c_\gamma Q(t)+c_\gamma\E_{\Ga
N}(t)+c_\gamma\sum_{k=2}^m\E^{k}(t),\quad 0\leq t\leq T,\ee and
\be\label{1.12n} Q(t)\leq c_\gamma Q(0)
+c_\gamma\int_0^t\e[\e(1+\E^{1/2}(t)\i)Q(t)+Q_{\Ga
N}(t)+\sum_{k=2}^m\E^k(t)\i]dt,\ee for $t\in[0,T]$, where $\E(t)$
and $Q(t)$ are given in $(\ref{1.28})$ and $(\ref{1.30})$,
respectively.
\end{thm}

{\bf Proof.}\,\,\,It will suffice to make some revisions on the
proofs of Lemmas 2.3 and 2.4, respectively.

Using the ellipticity that there is $c_\gamma>0$ such that
$$\|w\|^2_{k+1}\leq
c_\gamma\e(\|\B(t)w\|^2_{k-1}+\|w_{\nu_B}\|^2_{k-1/2,\,\Ga}+\|w\|^2_{k}\i),\quad
w\in H^k(\Om),$$ for $0\leq k\leq m-1$, in the proof of Lemma 2.3
yields the inequality (\ref{2.26n}).

Moreover, the second inequality (\ref{1.12n}) is based on the
following

\begin{lem}
Let $\gamma>0$ be given and $u$ be a solution of the problem
$(\ref{2.1n})$ on the interval $[0,T]$ for some $T>0$ such that the
condition $(\ref{1.25})$ holds true. Let $w\in
H^1\e((0,T)\times\Om\i)$ solve the linear
problem\be\label{1.4}\cases{\ddot{w}(t)=\B(t)w+F(t)\quad (t,x)\in
(0,T)\times\Om,\cr w_{\Ga_1}=0,\quad w_{\nu_B}|_{\Ga_0}=\var\quad
(t,x)\in(0,T)\times\Ga,\cr w(0)=w^0,\quad\dot{w}(0)=w^1,\quad
x\in\Om,}\ee where \be\label{x1} w_{\nu_B}=\<B\nabla w,\,\nu\>.\ee
Set$$ \Upsilon(t)=\|\dot{w}(t)\|^2+\|\nabla w(t)\|^2,\quad
\Upsilon_{\Ga N}(t)=\|\var(t)\|^2_{H^{1/2}(\Ga_0)}.$$ Then there is
$c_\gamma>0$ such that\be\label{1.5n} \Upsilon(t)\leq c_\gamma
\Upsilon(0)+c_\gamma\int_0^t\e[(1+\|\dot{u}(t)\|_{m-1})\Upsilon(\tau)+\Upsilon_{\Ga
N}(\tau)+\|F(\tau)\|^2\i]d\tau,\ee for $ 0\leq t\leq T$.
\end{lem}

{\bf Proof.}\,\,\,Let $$P(t)=\|\dot{w}\|^2+\e(B\nabla w,\nabla
w\i).$$ Then \be\label{1n}
\dot{P}(t)=2\e(F,\,\,\dot{w}\i)+\e(\dot{B}\nabla w,\nabla w\i)+
2\int_{\Ga_0}\var\dot{w}\Ga.\ee Using the estimate
$$|(\var,\dot{w})_{L^2(\Ga_0)}|\leq
\|\dot{w}\|_{H^{-1/2}(\Ga_0)}\|\var\|_{H^{1/2}(\Ga_0)}\leq
c\Upsilon(t)+c\Upsilon_{\Ga N}(t)$$ in (\ref{1n}) gives the
inequality (\ref{1.5n}).

\def\theequation{3.\arabic{equation}}
\setcounter{equation}{0}
\section{Locally exact controllability; the Dirichlet action }
\hskip\parindent The first step of the proof for the local exact
controllability depends on the following fact: Let $\X$ and $\Y$
be Banach spaces and $\Phi$: $\O\rw\Y$, where $\O$ is an open
subset of $\X$, be Frech$\acute{e}$t differentiable. If
$\Phi'(X_0)$: $\X\rw\Y$ is surjective, then there is an open
neighbourhood of $Y_0=\Phi(X_0)$ contained in the image
$\Phi(\O)$.

We start by specifying a value of $T$ about which we shall say more
later. We introduce a Banach space $\X^m_0(T)$ as follows.
$\X^m_0(T)$ consists of all
 the functions \be\label{3.2}
\var\in\cap_{k=0}^{m-2}C^k\e([0,T],H^{m-1/2-k}(\Ga_0)\i),\quad
\var^{(k)}\in H^1\e((0,T)\times\Ga_0)\i),\ee \be\label{3.3}
\var^{(k)}(0)=0,\quad x\in\Ga_0,\quad0\leq k\leq m-1,\ee with the
norm\be\label{3.4}
\|\var\|^2_{\X^m_0(T)}=\sum_{k=0}^{m-2}\|\var^{(k)}\|^2_{C\e([0,T],H^{m-k-3/2}(\Ga_0)\i)}
+\sum_{k=0}^{m-1}\|\var^{(k)}\|^2_{H^1\e((0,T)\times\Ga_0\i)}.\ee

Let an equilibrium solution $w\in H^m_{\Ga_1}(\Om)$ be given. We
invoke Theorem 1.1 to define a map for $\var\in\X^m_0(T)$ by
setting\be\label{3.7} \Phi(\var)=\e(u(T),\dot{u}(T)\i),\ee where $u$
is the solution of the following problem\be\label{3.8}
\cases{\ddot{u}=\sum a_{ij}(x,\nabla u)u_{x_ix_j}+b(x,\nabla u)\quad
(t,x)\in(0,T)\times\Om,\cr u|_{\Ga_1}=0,\quad t\in(0,T),\cr
u|_{\Ga_0}=w|_{\Ga_0}+\var,\quad t\in(0,T), \cr u(0)=w,\quad
\dot{u}(0)=0.}\ee Let $\varepsilon_T>0$ be given by Theorem 1.1.
Then\be\label{3.10} \Phi:\quad B_{\X^m_0(T)}(0,\varepsilon_T)\rw
\e(H^{m}(\Om)\cap H^1_{\Ga_1}(\Om)\i)\times \e(H^{m-1}(\Om)\cap
H^1_{\Ga_1}(\Om)\i),\ee where
$B_{\X^m_0(T)}(0,\varepsilon_T)\subset\X^m_0(T)$ is the ball with
the radius $\varepsilon_T$ centered at $0$. We observe that
$\Phi(0)=(w,0)$.

We need to evaluate\be\label{3.11} \Phi'(0)\var=\frac{\pl}{\pl
\si}\Phi(\si\var)|_{\si=0},\quad\var\in\X^m_0(T).\ee It is easy to
check that\be\label{3.12} \Phi'(0)\var=\e(v(T),\dot{v}(T)\i),\ee
where $v(t,x)$ is the solution of the linear system with variable
coefficients in the space variable \be\label{3.13}
\cases{\ddot{v}=\A v+F(v),\quad (t,x)\in(0,T)\times\Om,\cr
v|_{\Ga_1}=0,\quad t\in(0,T),\cr v|_{\Ga_0}=\var,\quad
t\in(0,T),\cr v(0)=\dot{v}(0)=0,}\ee where \be\label{3.14} \A
v=\sum_{ij=1}^na_{ij}(x,\nabla w)v_{x_ix_j},\ee \be\label{3.16}
F=(F_1,\cdots,F_n),\ee $$F_i=\sum_{lj} [a_{ljy_i}(x,\nabla
w)(w_{x_lx_j}-a_{ijy_l}(x,\nabla w)w_{x_lx_j}]+b_{y_i}(x,\nabla
w).$$

We now verify that $\Phi'(0)$ is surjective. In the language of
control theory surjection is just exact controllability, which for
a reversible system such as (\ref{3.13}) is equivalent to null
controllability.

Explicitly, one has to show that, for specified $T$, given $v^0\in
H^m(\Om)\cap H^1_{\Ga_1}(\Om)$ and $v^1\in H^{m-1}(\Om)\cap
H^1_{\Ga_1}(\Om)$,  one can find $\var\in\tilde{\X}^m_0(T)$ such
that the solution to \be\label{3.15} \cases{\ddot{v}=\A v+F(v),\quad
(t,x)\in(0,T)\times\Om,\cr v|_{\Ga_1}=0,\quad t\in(0,T),\cr
v|_{\Ga_0}=\var,\quad t\in(0,T),\cr v(0)=v^0,\quad\dot{v}(0)=v^1}\ee
satisfies \be\label{3.17}v(T)=\dot{v}(T)=0,\ee where
$\tilde{\X}^m_0(T)$ is the Banach space of all function with
(\ref{3.2}) and the norm $(\ref{3.4})$ but with  (\ref{3.3})
replaced by \be\label{3.3!}\var^{(k)}(T)=0,\quad
x\in\Ga_0,\quad0\leq k\leq m-1.\ee

Theorem 1.2 is then established by the following

\begin{thm}
Let an equilibrium $w\in H^m(\Om)\cap H^1_{\Ga_1}(\Om)$ be exactly
controllable. Let $T>T_0$ be given where $T_0$ is defined by
$(\ref{2.4*})$. Then, for any
$$(v^0,v^1)\in \e(H^m(\Om)\cap H^1_{\Ga_1}(\Om)\i)\times
\e(H^{m-1}(\Om)\cap H^1_{\Ga_1}(\Om)\i),$$ there is a
$\var\in\tilde{\X}^m_0(T)$ such that the solution
$$v\in\cap_{k=0}^mC^k\e([0,T], H^{m-k}(\Om)\i)$$ of the problem
$(\ref{3.15})$ satisfies $(\ref{3.17})$.
\end{thm}

{\bf Distributed Control.}\,\,\, As to the exact controllability of
linear systems by distributed control there is a long history and
the results are rich where many approaches are involved. Here the
distributed control means that solutions $(v(t),\dot{v}(t))$ of the
controlled system (\ref{3.15}) are only in the space $L^2(\Om)\times
H^{-1}(\Om)$ for $t\in[0,T]$. We just mention what we need in this
paper. One of the useful approaches is the multiplier method below,
introduced by Ho \cite{H} and Lions \cite{L}, to control the linear
system by its duality system.

We start with the wave equation \be\label{3.18}
\cases{\ddot{\phi}=\B\phi,\quad (t,x)\in(0,T)\times\Om,\cr
\phi|_\Ga=0,\quad t\in(0,T),\cr
\phi(0)=\phi_0,\quad\dot{\phi}(0)=\phi_1,}\ee where the operator
$\B$ is defined by\be\label{3.40} \B v=\A v-F(v)-v\, \div F,\quad
v\in H^2(\Om),\ee that is the dual system of the system
(\ref{3.15}). We have the following Green formula\be\label{3.42}
(v,\B u)=\e(\B^\star
v,\,u\i)+\int_\Ga\e[vu_{\nu_A}-uv_{\nu_A}-uv\<AF,\nu\>\i]d\Ga,\quad
u,\,v\in H^2(\Om),\ee where\be\label{3.41} \B^\star v=\A v+F(v).\ee

Given $(\phi_0,\phi_1)\in H^1_0(\Om)\times L^2(\Om)$, the problem
(\ref{3.18}) admits a unique solution. We then solve the problem
\be\label{3.19} \cases{\ddot{\psi}=\A \psi+F(\psi),\quad
(t,x)\in(0,T)\times\Om,\cr\psi(T)=\dot{\psi}(T)=0,\quad x\in\Om, \cr
\psi|_{\Ga_1}=0,\quad\psi|_{\Ga_0}=\phi_{\nu_A},\quad t\in(0,T),}\ee
where $\phi_{\nu_A}=\<A\nabla\phi,\nu\>$, $A=\e(a_{ij}(x,\nabla
w)\i)$, and $\phi$ is produced by $(\ref{3.18})$. Let $\psi$ be the
solution of the problem (\ref{3.19}). We then have constructed a
control $\phi_{\nu_A}$ on $(0,T)\times\Ga_0$ moving the initial
state $\e(\psi(0),\dot{\psi}(0)\i)$ to rest at the time $T$.

We define a mapping $\Lam$: $ H^1_0(\Om)\times L^2(\Om)\rw
H^{-1}(\Om)\times L^2(\Om)$ by \be\label{3.25} \Lam(\phi_0,
\phi_1)=(\dot{\psi}(0), -\psi(0)).\ee A formal use of Green's
formula yields, after we multiply (\ref{3.18}) by $\phi$ and
integrate by parts over $\Sigma=(0,T)\times\Om$, \be\label{3.20}
\<\Lam(\phi_0, \phi_1), (\phi_0, \phi_1)\>_{L^2(\Om)\times
L^2(\Om)}=\int_{\wp_0}\phi_{\nu_A}^2\wp,\ee where
$\wp_0=(0,T)\times\Ga_0$.

Let constants $c_1>0$ and $c_2>0$ be such that\be\label{3.21}
c_1E(\phi_0,\phi_1)\leq \int_{\wp_0}\phi_{\nu_A}^2d\wp\leq
c_2E(\phi_0,\phi_1),\ee for $(\phi_0,\phi_1)\in H^1_0(\Om)\times
L^2(\Om)$ where\be\label{3.31}
E(\phi_0,\phi_1)=\|A^{1/2}\nabla\phi_0\|^2+\|\phi_1\|^2.\ee Then,
for any $(\phi_0,\phi_1)\in H^1_0(\Om)\times H^1(\Om)$, one has
$\phi_{\nu_A}\in L^2\e((0,T)\times\Ga_0\i)$ that drives the system
starting from $(\psi(0),\dot{\psi}(0))$ at the time $t=0$ to rest at
the time $T$.

Then the key point is to establish the inequality (\ref{3.21}). For
$\A$ being the classical Laplacian and $F=0$, the inequality
(\ref{3.21}) was proved in Ho \cite{H}. For $\A$ with variable
coefficients in space, such as (\ref{3.14}), and $F=0$, the
inequality (\ref{3.21}) was established under some geometric
conditions in Yao \cite{Y}, where the geometrical method was
introduced. Without geometric conditions, the inequality
(\ref{3.21}) is not true even if the control portion $\Ga_0$ of
$\Ga$ is the whole boundary. A counterexample was given by Yao
\cite{Y}. Then the geometrical method was extended by Lasiecka,
Triggiani, and Yao \cite{LTY} to include the case of the first order
terms $F\not=0$. This method was again extended to study the
modeling and control problems of thin shells by Chai etc.,
\cite{C1}, \cite{C2}, and Lasiecka, etc., \cite{LT1}. A recent
survey paper on the geometrical method is by Gulliver, etc.,
\cite{GU}.

The lemma below follows by Lasiecka, Triggiani, and Yao \cite{LTY},
Theorem 3.2, where a uniqueness result, needed, is provided by
Triggiani and Yao \cite{TY}, Theorem 10.1.1.

\begin{lem}
Let an equilibrium $w\in H^m(\Om)\cap H^1_{\Ga_1}(\Om)$ be exactly
controllable and $T_0$ be given by the formula $(\ref{2.4*})$. Then,
for $T>T_0$ given, $\Lam$ is an isomorphism from $H^1_0(\Om)\times
L^2(\Om)$ onto $H^{-1}(\Om)\times L^2(\Om)$. In particular, there
are $c_1>0$ and $c_2>0$ such that the inequality $(\ref{3.21})$
holds true.
\end{lem}

However, the above control strategy only gives distributed control
functions because solutions $(\psi(t),\dot{\psi}(t))$ of the
controlled system (\ref{3.19}) are only in $L^2(\Om)\times
H^{-1}(\Om)$ no matter $(\phi_0,\phi_1)$ are smooth or not.
Indeed, since $\phi_{\nu_A}(T)\not=0$ for any $x\in\Ga_0$, the
compatible
condition $\psi(T)=\phi_{\nu}(T)$ for $x\in\Ga_0$ is never true.\\

{\bf Smooth Control.}\,\,\,We shall modify the above control
strategy to obtain smooth controls to meet the need of Theorem
3.1.

Let $k\geq 1$ be an integer. Let $\Xi_0^k(\Om)$ consist of the
functions $u$ in $H^k(\Om)$ with the boundary
conditions\be\label{3.24}\cases{\B^iu|_{\Ga_0}=0,\quad 0\leq i\leq
l-1,\quad\mbox{if}\quad k=2l;\cr \B^iu|_{\Ga_0}=0,\quad 0\leq
i\leq l\quad\mbox{if}\quad k=2l+1,}\ee and with the norms of
$H^k(\Om)$ where $\B$ is given by (\ref{3.40}).

Let $T_0$ be given by the formula (\ref{2.4*}) and $T>T_1>T_0$ be
given. We assume that $z\in C^{\infty}(-\infty,\infty)$ is such
that $0\leq z(t)\leq 1$ with \be\label{3.22} z(t)=\cases{0,\quad
t\geq T,\cr 1,\quad t\leq T_1.}\ee

For $(\phi_0,\phi_1)\in \Xi_0^{m+1}(\Om)\times\Xi_0^{m}(\Om)$
given, we solve the problem (\ref{3.18}) and then, in stead of
(\ref{3.19}), we solve the following problem \be\label{3.23}
\cases{\ddot{\psi}=\A \psi+F(\psi),\quad
(t,x)\in(0,T)\times\Om,\cr\psi(T)=\dot{\psi}(T)=0,\quad x\in\Om,
\cr \psi|_{\Ga_1}=0,\quad\psi|_{\Ga_0}=z\phi_{\nu_A},\quad
t\in(0,T).}\ee

Let $\Lam$ be given by (\ref{3.25})  where $\psi$ in $(\ref{3.25})$
are solutions of the problem $(\ref{3.23})$ this time. It is easy to
check that, for any $(\phi_0,\phi_1)$, $(\var_0,\var_1)\in
H^1_0(\Om)\times L^2(\Om)$,\be\label{3.26}
\<\Lam(\phi_0,\phi_1),(\var_0,\var_1)\>_{L^2(\Om)\times
L^2(\Om)}=\int_{\wp_0}z(t)\phi_{\nu_A}\var_{\nu_A}d\wp,\ee with $
\wp_0=(0,T)\times\Ga_0$, where $\phi$ and $\var$ are solutions of
the problem (\ref{3.18}) with initial data $(\phi_0,\phi_1)$ and
$(\var_0,\var_1)$, respectively.

We shall show that the problem (\ref{3.23}) provides smooth
controls to Theorem 3.1 by the following lemma.

\begin{lem}
Let $k\geq0$ be an integer and $\Lam$ be given by $(\ref{3.25})$
where $\psi$ is the solution of the problem $(\ref{3.23})$. There
are then $c_1>0$ and $c_2>0$ such that$$
c_1\|(\phi_0,\phi_1)\|_{H^{k+1}(\Om)\times
H^k(\Om)}\leq\|\Lam(\phi_0,\phi_1)\|_{H^{k-1}(\Om)\times
H^{k}(\Om)}\leq c_2\|(\phi_0,\phi_1)\|_{H^{k+1}(\Om)\times
H^k(\Om)},$$ \be\label{3.27}\any\,\,
(\phi_0,\phi_1)\in\Xi_0^{k+1}(\Om)\times\Xi_0^k(\Om).\ee In
particular, $\Lam$  are  isomorphisms from
$\Xi_0^{2}(\Om)\times\Xi_0^1(\Om)$ onto $L^2(\Om)\times
H^1_{\Ga_1}(\Om)$ and from $\Xi_0^{k+1}(\Om)\times\Xi_0^k(\Om)$ onto
$\e(H^{k-1}(\Om)\cap H^1_{\Ga_1}(\Om)\i)\times \e( H^k(\Om)\cap
H^1_{\Ga_1}(\Om)\i)$ for $k\geq2$, respectively.
\end{lem}

{\bf Proof.}\,\,\,Lemma 3.1 shows that the inequality (\ref{3.27})
is true for $k=0$.

We now proceed to prove the inequality (\ref{3.27}) by induction
on $k$. Let the inequality (\ref{3.27}) be true for some integer $
k\geq 0$. We want to show that the inequality (\ref{3.27}) hold
with $k$ replaced by $k+1$.

{ \bf Case I}\,\,\,Let $k=2l$ for some $l\geq1$.

Let
\be\label{3.36}(\phi_0,\phi_1)\in\Xi_0^{k+2}(\Om)\times\Xi_0^{k+1}(\Om)\ee
be given. Suppose that $\phi$ is the solution of the problem
(\ref{3.18}) corresponding to the initial data $(\phi_0,\phi_1)$.
Then $\phi^{(2i)}$ and $\phi^{(2i+1)}$ are the solutions of the
problem (\ref{3.18}) corresponding to the initial data
$(\B^i\phi_0,\B^i\phi_1)$ and $(\B^i\phi_1,\B^{i+1}\phi_0)$,
respectively, for $0\leq i\leq l$, where $\B$ is given by
(\ref{3.40}).

For any
$(\var_0,\var_1)\in\Xi_0^{2(k+1)}(\Om)\times\Xi_0^{2k+1}(\Om)$, let
$\var$ be the solution of the problrm (\ref{3.18}) with the initial
data $(\var_0,\var_1)$. Then $\var^{(2i)}$ and $\var^{(2i+1)}$ are
the solutions of the problem (\ref{3.18}) corresponding to the
initial data $(\B^i\var_0,\B^i\var_1)$ and
$(\B^i\var_1,\B^{i+1}\var_0)$, respectively, for $0\leq i\leq k$.
Using the initial data $(\phi_0,\phi_1)$ and
$(\B^{k+1}\var_0,\B^{k+1}\var_1)$ in the formula (\ref{3.26}), we
obtain \be\label{3.43}
\e(\psi(0),\B^{k+1}\var_1\i)-\e(\dot{\psi}(0),\,\B^{k+1}\var_0\i)
=-\int_{\wp_0}z(t)\phi_{\nu_A}\var^{(2k+2)}_{\nu_A}d\wp.\ee

In one hand, by integration by parts with respect to the variable
$t$ on $[0,T]$, we obtain\beq\label{3.44}
&&-\int_{\wp_0}z(t)\phi_{\nu_A}\var^{(2k+2)}_{\nu_A}d\wp\nonumber\\
&&=\sum_{j=1}^{k}(-1)^j\phi^{(j)}_{\nu_A}(0)\var^{(2k+1-j)}_{\nu_A}(0)
+ \int_{\wp_0}\e(z(t)\phi_{\nu_A}\i)^{(k+1)}
\var^{(k+1)}_{\nu_A}d\wp\nonumber\\
&&=\sum_{j=0}^{l}\e(\B^{j}\phi_0\i)_{\nu_A}\e(\B^{k-j}\var_1\i)_{\nu_A}
-\sum_{j=0}^{l-1}\e(\B^{j}\phi_1\i)_{\nu_A}\e(\B^{k-j}\var_0\i)_{\nu_A}+I(\phi,\var),\eeq
where\be\label{3.46}I(\phi,\var)=\sum_{j=1}^{k+1}
\int_{\wp_0}z^{(j)}(t)\phi_{\nu_A}^{(k+1-j)}
\var^{(k+1)}_{\nu_A}d\wp+
\int_{\wp_0}z(t)\phi_{\nu_A}^{(k+1)}\var^{(k+1)}_{\nu_A}d\wp.\ee

On the other hand, using the formula (\ref{3.42}), the boundary
conditions (\ref{3.24}), and the equation (\ref{3.23}), we
obtain\beq\label{3.48}\e(\psi(0),\B^{k+1}\var_1\i)
&&=\e((\B^\star)^l\psi(0),\B^{l+1}\var_1\i)+
\sum_{j=0}^{l-1}\int_{\Ga}(\B^\star)^{j}\psi(0)\e(\B^{k-j}\var_1\i)_{\nu_A}d\Ga\nonumber\\
&&=-\e(A\nabla(\B^\star)^l\psi(0),\nabla\B^{l}\var_1\i)+
\sum_{j=0}^{l}\int_{\Ga_0}\psi^{(2j)}(0)\e(\B^{k-j}\var_1\i)_{\nu_A}d\Ga\nonumber\\
&&\quad-\e((\B^\star)^l\psi(0),F(\B^l\var_1)+(\B^l\var_1)\div
F\i),\eeq and\beq\label{3.49}
\e(\dot{\psi}(0),\B^{k+1}\var_0\i)&&=\e((\B^\star)^l\dot{\psi}(0),\B^{l+1}\var_0\i)+\sum_{j=0}^{l-1}
\int_{\Ga_0}\psi^{(2j+1)}(0)\e(\B^{k-j}\var_0\i)_{\nu_A}d\Ga.\eeq

Noting that $\psi^{(2j)}(0)=\phi^{(2j)}(0)=\B^j\phi_0$ and
$\psi^{(2j+1)}(0)=\B^j\phi_1$ on $\Ga_0$ and using
(\ref{3.43})-(\ref{3.49}), we have the following identity
\beq\label{3.54}&&
-\e(A\nabla(\B^\star)^l\psi(0),\nabla\B^{l}\var_1\i)-
\e((\B^\star)^l\dot{\psi}(0),\B^{l+1}\var_0\i)\nonumber\\
&&=I(\phi,\var)+\e((\B^\star)^l\psi(0),\,\,F(\B^l\var_1)+(\B^l\var_1)\div
F\i)\eeq

Since $\Xi_0^{2(k+1)}(\Om)\times\Xi_0^{2k+1}(\Om)$ is dense in
$\Xi_0^{k+2}(\Om)\times\Xi_0^{k+1}(\Om)$, the identity
(\ref{3.54}) is actually true for all
$(\var_0,\var_1)\in\Xi_0^{k+2}(\Om)\times\Xi_0^{k+1}(\Om)$.

Letting $\var_0=0$ in (\ref{3.54}), we obtain \be\label{3.55}
\e(\A(\B^\star)^l\psi(0),\B^l\var_1\i)=I(\phi,\var)+\e((\B^\star)^l\psi(0),\,\,F(\B^l\var_1)+(\B^l\var_1)\div
F\i),\ee for $\var_1\in\Xi_0^{k+1}(\Om)$ where $\var$ is the
solution of the problem (\ref{3.18}) for the initial data
$(0,\var_1)$. It is easy to check by the maximum principle for the
elliptic operator that \be\label{3.56} \ol{Image(\B^l)}=L^2(\Om).\ee

Moreover, by virtue of the inequality (\ref{3.21}) and Lemma 3.1, we
have the
estimate\beq\label{3.45} &&|I(\phi,\var)|\nonumber\\
&&\leq
c\sum_{j=0}^l\e[\int_0^T\int_{\Ga_0}\e((\phi^{(2j)}_{\nu_A})^2+(\phi^{(2j+1)}_{\nu_A})^2\i)d\Ga
dt\i]^{1/2}\e(\int_{\wp_0}(\var^{(k+1)}_{\nu_A})^2d\wp\i)^{1/2}\nonumber\\
&&\leq
c\sum_{j=0}^l\e(E(\B^j\phi_0,\B^j\phi_1)+E(\B^j\phi_1,\B^{j+1}\phi_0)\i)^{1/2}\|A^{1/2}\nabla\B^l\var_1\|\nonumber\\
&&\leq
c\e(\|\phi_0\|^2_{k+2}+\|\phi_1\|^2_{k+1}\i)^{1/2}\|\B^l\var_1\|_1.\eeq
In terms of (\ref{3.55})-(\ref{3.45}), we obtain
\beq\label{3.57}\|\A(\B^\star)^l\psi(0)\|_{-1}&&=\sup_{\|\B\var_1\|_1=1}\e(\A(\B^\star)^l\psi(0),\,\B^l\var_1\i)\nonumber\\
&&\leq
c\e(\|\phi_0\|^2_{k+2}+\|\phi_1\|^2_{k+1}\i)^{1/2}+c\|\psi(0)\|_k.\eeq

Furthermore, on the boundary $\Ga$ the problem (\ref{3.23})
implies\beq\label{3.58}
\|(\B^\star)^i\psi(0)\|_{H^{k+1/2-2i}(\Ga)}&&=\|\psi^{(2i)}(0)\|_{H^{k+1/2-2i}(\Ga)}
=\|\phi^{(2i)}_{\nu_A}(0)\|_{H^{k+1/2-2i}(\Ga_0)}\nonumber\\
&&=\|(\B^i\phi_0)_{\nu_A}\|_{H^{k+1/2-2i}(\Ga_0)}\leq
c\|\phi_0\|_{k+2},\quad 0\leq i\leq l.\eeq

Now, using the ellipticity of the operator $\B^\star$ and from
(\ref{3.57}) and (\ref{3.58}), we have\beq\label{3.56*}
\|\psi(0)\|_{k+1}&&\leq
c\|\B^\star\psi\|_{k-1}+c\|\psi(0)\|_{H^{k+1/2}(\Ga)}+c\|\psi(0)\|_k\nonumber\\
&&\leq
c\|(\B^\star)^2\psi\|_{k-3}+c\|\B^\star\psi(0)\|_{H^{k+1/2-2}(\Ga_0)}+\|\phi_0\|_{k+2}+c\|\psi(0)\|_k\nonumber\\
&&\leq
c\|(\B^\star)^l\psi(0)\|_{-1}+c\|\phi_0\|_{k+2}+c\|\psi(0)\|_k\nonumber\\
&&\leq c\e(\|\phi_0\|^2_{k+2}+\|\phi_1\|^2_{k+1}\i)^{1/2},\eeq where
the induction assumption $\|\psi(0)\|_k\leq c
\e(\|\phi_0\|^2_{k+1}+\|\phi_1\|^2_{k}\i)^{1/2}$ is used.

A similar argument yields\be\label{3.59} \|\dot{\psi}(0)\|_k\leq c
\e(\|\phi_0\|^2_{k+2}+\|\phi_1\|^2_{k+1}\i)^{1/2},\ee after we let
$\var_0\in\Xi_0^{k+2}(\Om)$ and $\var=0$ in (\ref{3.54}).

Next, let us prove the left hand side of the inequality (\ref{3.27})
where $k$ is replaced by $k+1$. We set $\var_0=\phi_0$ and
$\var_1=\phi_1$ in (\ref{3.54}) and use Lemma 3.1  to
obtain\beq\label{3.60}
&&c\e(\|\psi(0)\|^2_{k+1}+\|\dot{\psi}(0)\|^2_k\i)^{1/2}
E^{1/2}(\B^l\phi_1,\B^{l+1}\phi_0)\nonumber\\
&&\geq
I(\phi,\phi)-c_1\|\psi(0)\|_k\|A^{1/2}\nabla\B^l\phi_1\|\nonumber\\
&&\geq \int_0^{T_1}\int_{\Ga_0}\e(\phi^{(k+1)}_{\nu_A}\i)^2d\Ga
dt-\varepsilon\int_{\wp_0}\e(\phi^{(k+1)}_{\nu_A}\i)^2d\wp-
c_\varepsilon\sum_{j=0}^{k}\int_{\wp_0}\e(\phi^{(j)}_{\nu_A}\i)^2d\wp\nonumber\\
&&\quad-\varepsilon\|\psi(0)\|^2_{k}-c_\varepsilon\|A^{1/2}\nabla\B^l\phi_1\|^2\nonumber\\
&&\geq
c_1E(\B^l\phi_1,\B^{l+1}\phi_0)-c_2\e(\|\phi_0\|^2_{k+1}+\|\phi_1\|^2_k\i).\eeq
In addition,
$(\phi_0,\phi_1)\in\Xi_0^{k+2}(\Om)\times\Xi_0^{k+1}(\Om)$ implies,
 by the ellipticity of the operator $\B$,\be\label{3.61}
\|\phi_0\|^2_{k+2}+\|\phi_1\|^2_{k+1}\leq
cE(\B^l\phi_1,\B^{l+1}\phi)+c\e(\|\phi_0\|^2_{k+1}+\|\phi_1\|^2_k\i).\ee
Then the inequalities (\ref{3.60}) and (\ref{3.61}) give, via the
induction assumption $\|\psi(0)\|_k^2+\|\dot{\psi}(0)\|^2_{k-1}\geq
c\e(\|\phi_0\|^2_{k+1}+\|\phi_1\|^2_{k}\i)$, \be\label{3.62}
\|\psi(0)\|^2_{k+1}+\|\dot{\psi}(0)\|^2_k\geq
c\e(\|\phi_0\|^2_{k+2}+\|\phi_1\|^2_{k+1}\i).\ee

The relations (\ref{3.56}), (\ref{3.59}) and (\ref{3.62}) mean that
the inequality (\ref{3.27}) is true with $k$ replaced by $k+1$ if
$k=2l$ for some $l\geq1$.

{\bf Case II}\,\,\,If $k=2l+1$, a similar argument can establish the
inequality (\ref{3.27}) where $k$ is replaced by $k=1$.

Then Lemma 3.2 follows by induction.

\begin{lem}
Let $\phi$ solve the problem $(\ref{3.18})$ with the initial data
$(\phi_0,\phi_1)\in\Xi_0^2(\Om)\times\Xi_0^1(\Om)$.
Then\be\label{3.65} \phi_{\nu_A}\in C\e([0,T],H^{1/2}(\Ga)\i)\cap
H^1\e((0,T)\times\Ga\i).\ee
\end{lem}

{\bf Proof.}\,\,\,For any $T>0$ given, there is $c_T>0$ such that
$$\|\phi(t)\|_2\leq c_T\e(\|\phi_0\|_2^2+\|\phi_1\|^2_1\i)\quad\any\,t\in[0,T],$$ which
implies $\phi_{\nu_A}\in C\e([0,T],H^{1/2}(\Ga)\i)$.

Since $\phi'$ is the solution of the problem (\ref{3.18}) for the
initial data $(\phi_1,\B\phi_0)\in H^1_0(\Om)\times L^2(\Om)$, Lemma
3.1 implies $\phi'_{\nu_A}\in L^2\e((0,T)\times\Ga\i)$.

To complete the proof, it is remaining to show that $\phi_{\nu_A}\in
L^2\e((0,T),H^1(\Ga)\i)$.

Let $X$ be a vector field of the manifold $\Ga$, that is,
$X(x)\in\Ga_x$ for each $x\in\Ga$. We extend $X$ to the whole
$\ol{\Om}$ to be a vector field on the manifold $(\ol{\Om},g)$ where
$g=A^{-1}=\e(a_{ij}(x,\nabla w)\i)^{-1}$.

Let\be\label{3.66} v=X(\phi),\quad (t,x)\in(0,T)\times\Om.\ee Then
$v$ solves the problem\be\label{3.67}\cases{\ddot{v}=\B
v+[X,\B]\phi,\quad (0,T)\times\Om,\cr v|_\Ga=0,\quad t\in(0,T), \cr
v(0)=X(\phi_0)\in H^1_0(\Om),\quad \dot{v}(0)=X(\phi_1)\in
L^2(\Om),}\ee where $[X,\B]\phi=X(\B\phi)-\B X(\phi)$ with the
estimate \be\label{3.71}\|[X,\B]\phi(t)\|\leq
c_T\e(\|\phi_0\|_2^2+\|\phi_1\|^2_1\i)\quad\any\,t\in[0,T].\ee

Let $H$ be a vector field on $\ol{\Om}$ with $$H(x)=\nu_A\quad
x\in\Ga.$$ We multiply the both sides of the equation in
(\ref{3.67}) by $H(v)$ and integrate by parts over
$\Sigma=(0,T)\times\Om$ to obtain\beq\label{3.68}
&&\int_\wp\e[v_{\nu_A}^2-\frac{1}{2}|\nabla_{g}v|^2_{g}|\nu_A|^2_{g}\i]d\wp\nonumber\\
&&=\e(\dot{v},H(v)\i)|_0^T+\int_\Sigma\e[D_{g}H(\nabla_{g}v,\nabla_{g}v)+\frac{1}{2}\e(\dot{v}^2
-|\nabla_{g}v|_{g}^2\i)\div H\i]d\Sigma\nonumber\\
&&\quad+\e(F(v)+v\div F-[X,\B]\phi,\,\,\,H(v)\i),\eeq where
$\wp=(0,T)\times\Ga$. In addition, the boundary condition $v|_\Ga=0$
implies
\be\label{3.69}|\nabla_{g}v|^2_{g}=\frac{1}{|\nu_A|^2_{g}}v_{\nu_A}^2\quad\any\,x\in\Ga.\ee
In terms of (\ref{3.68}), (\ref{3.69}) and (\ref{3.71}), we
obtain\be\label{3.70} \int_\wp v^2_{\nu_A}d\wp\leq
c_T\e(\|\phi_0\|^2_2+\|\phi_1\|^2_1\i).\ee

Since
$$v_{\nu_A}=\nu_A\e(X(\phi)\i)=X(\phi_{\nu_A})+[\nu_A,X]\phi\quad
x\in\Ga,$$ by (\ref{3.70}), we have
$$\int_0^T\int_\Ga|X(\phi_{\nu_A})|^2d\Ga dt\leq
c_{T,X}\e(\|\phi_0\|^2_2+\|\phi_1\|^2_1\i),$$ for any vector field
$X$ of the manifold $\Ga$, that is, $\phi_{\nu_A}\in
L^2\e((0,T),H^1(\Ga)\i)$.  {\bf\textbf{$\|$}}\\

{\bf The Proof of Theorem 3.1}\,\,\,Let $(v_0,v_1)\in
\e(H^m(\Om)\cap H^1_{\Ga_1}(\Om)\i)\times \e(H^{m-1}(\Om)\cap
H^1_{\Ga_1}(\Om)\i)$ be given. By Lemma 3.2, there is
$(\phi_0,\phi_1)\in\Xi_0^{m+1}(\Om)\times\Xi_0^m(\Om)$ such that the
control $\var=z\phi_{\nu_A}$ on $\wp_0=(0,T)\times\Ga_0$ drives the
system (\ref{3.15}) to rest at the time $T$, where $\phi$ is the
solution of the problem (\ref{3.18}) with the initial data
$(\phi_0,\phi_1)$.

Since $\phi^{(k)}$ are the solutions of the problem (\ref{3.18})
with the initial
data\be\label{3.64}\cases{(\B^l\phi_0,\B^l\phi_1)\quad\mbox{if}\quad
k=2l\,\,\mbox{or},\cr (\B^l\phi_1,\B^{l+1}\phi_0)\quad\mbox{if}\quad
k=2l+1,}\ee for $0\leq k\leq m-1$, Lemma 3.3 implies
$\var=z\phi_{\nu_A}\in\tilde{\X}^m_0(\Om)$.

\def\theequation{4.\arabic{equation}}
\setcounter{equation}{0}
\section{Locally exact controllability; the Neumann action }
\hskip\parindent  Let $T>0$ be given. This time, we introduce a
Banach space $\X^m_{0 N}(T)$ as follows. $\X^m_{0 N}(T)$ consists of
all the functions \be\label{3.2n}
\var\in\cap_{k=0}^{m-2}C^k\e([0,T],H^{m-k-3/2}(\Ga_0)\i),\quad
\var^{(k)}\in L^2\e((0,T),H^{1/2}(\Ga_0)\i),\ee \be\label{3.3n}
\var^{(k)}(0)=0,\quad x\in\Ga_0,\quad0\leq k\leq m-1,\ee with the
norm\be\label{3.4n} \|\var\|^2_{\X^m_{0
N}(T)}=\sum_{k=0}^{m-2}\|\var^{(k)}\|^2_{C\e([0,T],H^{m-1/2-k}(\Ga_0)\i)}
+\sum_{k=0}^{m-1}\|\var^{(k)}\|^2_{L^2\e((0,T),
H^{1/2}(\Ga_0)\i)}.\ee

Let $w\in H^{m+1}(\Om)\cap H^1_{\Ga_1}(\Om)$ be given. We invoke
Theorem 1.4 to define a map for $\var\in\X^m_{0 N}(T)$ by
setting\be\label{3.7n} \Phi_N(\var)=\e(u(T),\dot{u}(T)\i),\ee where
$u$ is the solution of the following problem\be\label{3.8n}
\cases{\ddot{u}=\div \textbf{a}(x,\nabla u)\quad
(t,x)\in(0,T)\times\Om,\cr u|_{\Ga_1}=0,\quad t\in(0,T),\cr
\<\textbf{a}(x,\nabla u),\,\nu\>|_{\Ga_0}=\<\textbf{a}(x,\nabla
w),\,\nu\>+\var,\quad t\in(0,T), \cr u(0)=w,\quad \dot{u}(0)=0.}\ee
Let $\varepsilon_T>0$ be given by Theorem 1.4. Then\be\label{3.10n}
\Phi_N:\quad B_{\X^m_{0 N}(T)}(0,\varepsilon_T)\rw \e(H^{m}(\Om)\cap
H^1_{\Ga_1}(\Om)\i)\times \e(H^{m-1}(\Om)\cap
H^1_{\Ga_1}(\Om)\i),\ee where $B_{\X^m_{0
N}(T)}(0,\varepsilon_T)\subset\X^m_{0 N}(T)$ is the ball with the
radius $\varepsilon_T$ centered at $0$.

We observe that, since $w\in H^{m+1}(\Om)\cap H^1_{\Ga_1}(\Om)$,
$$\Phi_N(0)=(w,0)\in \e(H^{m+1}(\Om)\cap
H^1_{\Ga_1}(\Om)\i)\times \e(H^{m}(\Om)\cap H^1_{\Ga_1}(\Om)\i).$$
Then Theorem 1.5 is equivalent to the following claim: For some
$T>0$ there are $\varepsilon_1>0$ and $\varepsilon_2>0$ with
$\varepsilon_T\geq\varepsilon_2$ such that \be\label{3.10nn}
B_{H^{m+1}(\Om)\times
H^{m}(\Om)}\e((w,0),\varepsilon_1\i)\subset\Phi_N\e(B_{\X^m_{0
N}(T)}(0,\varepsilon_2)\i),\ee where $B_{H^{m+1}(\Om)\times
H^{m}(\Om)}\e((w,0),\varepsilon_1\i)$ is the ball with the radius
$\varepsilon_1$ centered at $0$ in the space $H^{m+1}(\Om)\times
H^{m}(\Om)$.

The map $\Phi_N$ is $Fr\acute{e}chet$ differentiable on $B_{\X^m_{0
N}(T)}(0,\varepsilon_T)$. In particular, \be\label{3.11n}
\Phi_N'(0)\var=\e(v(T),\dot{v}(T)\i),\quad\var\in\X^m_{0 N}(T),\ee
where $v(t,x)$ is the solution of the linear system with variable
coefficients in the space variable \be\label{3.13n}
\cases{\ddot{v}=\div A(x,\nabla w)\nabla v,\quad
(t,x)\in(0,T)\times\Om,\cr v|_{\Ga_1}=0,\quad t\in(0,T),\cr
v_{\nu_A}|_{\Ga_0}=\var,\quad t\in(0,T),\cr v(0)=\dot{v}(0)=0,}\ee
where $v_{\nu_A}=\<A(x,\nabla w)\nabla v,\,\nu\>$.

The proof of the exact controllability with the Neumann action
depends on the fact:

\begin{pro}\label{p1}
Let $\X_1$, $\X_2$,  $\Y_1$, and $\Y_2$ be Banach spaces with
$\X_2\subset\X_1$, $\Y_2\subset\Y_1$,  $\ol{\X_2}=\X_1$, and
$\ol{\Y_2}=\Y_1$.  Suppose that $\Phi:$ $\ball_{\X_i}(0,r)\rw
\Y_i$ are mappings of $C^1$ for $i=1$, $2$ such that
\be\label{nc1} \Y_2\subset\Phi'(0)\X_1.\ee There is
$\varepsilon>0$ such that\be\label{nc2}
\ball_{\Y_2}\e(\Phi(0),\varepsilon\i)\subset\Phi\e(\ball_{\X_1}(0,r)\i).\ee
\end{pro}

{\bf Proof.}\,\,\,Let $y_0=\Phi(0)$. It will suffice to prove that
for any $y$ in $\Y_2$ near $y_0$, the equation\be\label{nc3}
\Phi(x)-\Phi(0)=y-y_0\ee has a solution $x$ in
$\ball_{\X_1}(0,r)$. This can be done by a modification of the
proof of Theorem (3.1.19) in Berger \cite{B}.

We denote by $\X_1/\ker\Phi'(0)$ the quotient space where
$$\ker\Phi'(0)=\{\,x\,|\,x\in\X_1,\Phi'(0)x=0\,\}.$$ The assumptions (\ref{nc1})
imply that the inversion of $\Phi'(0)$: $\Y_2\rw\X_1/\ker\Phi'(0)$
exists, is closed, and therefore is bounded. Then, there is $C>0$
such that\be\label{nc4} C\|\Phi'(0)x\|_{\Y_2}\geq
d\e(x,\ker\Phi'(0)\i)\ee for $x\in\X_1$ such that
$\Phi'(0)x\in\Y_2$ where $d\e(x,\ker\Phi'(0)\i)$ is the distance
from $x$ to the space $\ker\Phi'(0)$ in $\X_1$.

Now we can construct a sequence $\{\,x_k\,\}$ as follows. Let
$\varepsilon>0$ be given. Let\be\label{nc5}
R(x)=\Phi(x)-\Phi(0)-\Phi'(0)x,\quad x\in\ball_\X(0,r).\ee Since
$\Phi$: $\ball_{\X_2}(0,r)\rw\Y_2$ is $C^1$ and $\ol{\X_2}=\X_1$,
we take $x_0\in\ball_{\X_1}(0,r)\cap\X_2$. Then $R(x_0)\in\Y_2$.
Next, the relations (\ref{nc1}) and (\ref{nc4}) imply that there
is $x*_1\in\X_1$ such
that\be\label{nc6}\Phi'(0)x*_1=y-y_0-R(x_0),\ee
\be\label{nc7}\|x*_1\|_{\X_1}\leq C\|y-y_0-R(x_0)\|_{\Y_2}.\ee If
$x*_1=x_0$, then $x_0$ is a solution to the equation (\ref{nc3})
and the constructing ends. We assume that $x*_1\not=x_0$. We take
$x_1\in\X_2$ such that\be\label{nc9}
\|x_1-x*_1\|_{\X_1}\leq\varepsilon\|x_1-x_0\|_{\X_1},\ee
\be\label{nc10} \|x_1\|_{\X_1}\leq C\|y-y_0-R(x_0)\|_{\Y_2}.\ee
Proceeding this procedure, we obtain two sequences
$\{\,x*_k\,\}\subset\X_1$ and $\{\,x_k\,\}\subset\X_2$
satisfying\be\label{nc11} \Phi'(0)x*_k=y-y_0-R(x_{k-1}),\ee
\be\label{nc12} \|x_k\|_{\X_1}\leq
C\|y-y_0-R(x_{k-1})\|_{\Y_2},\ee \be\label{nc13}
\|x_k-x*_k\|_{\X_1}\leq\varepsilon\|x_k-x_{k-1}\|_{\X_1},\ee for
$k\geq 1$.

A similar argument as in the proof of Theorem (3.1.19) in Berger
\cite{B} completes the proof. {\bf\textbf{$\|$}}\\

For $i=1$, $2$, let $$\X_i=\X^{m+i-1}_{0 N}(T),$$
$$\Y_i=\e(H^{m+i-1}(\Om)\cap
H^1_{\Ga_1}(\Om)\i)\times \e(H^{m+i-2}(\Om)\cap
H^1_{\Ga_1}(\Om)\i).$$ It is easy to check by Theorem \ref{t2.2}
that the mappings $\Phi_N$, given by (\ref{3.7n}), are of $C^1$
from $\ball_{\X_i}(0,r)$ to $\Y_i$ for $i=1$, $2$, and for some
$r>0$. By Proposition \ref{p1}, to prove Theorem 1.5 is to
establish the exact controllability of the system (\ref{3.13n}) on
the space $\e(H^{m+1}(\Om)\cap H^1_{\Ga_1}(\Om)\i)\times
\e(H^{m}(\Om)\cap H^1_{\Ga_1}(\Om)\i)$, which for a reversible
system such as (\ref{3.13n}) is equivalent to null
controllability.

Explicitly, one has to show that, for specified $T$, given
$(v_0,v_1)\in \e(H^{m+1}(\Om)\cap H^1_{\Ga_1}(\Om)\i)\times
\e(H^{m}(\Om)\cap H^1_{\Ga_1}(\Om)\i)$, one can find
$\var\in\tilde{\X}^m_{0 N}(T)$ such that the solution to
\be\label{3.15n} \cases{\ddot{v}=\div A(x,\nabla w)\nabla v,\quad
(t,x)\in(0,T)\times\Om,\cr v|_{\Ga_1}=0,\quad t\in(0,T),\cr
v_{\nu_A}|_{\Ga_0}=\var,\quad t\in(0,T),\cr
v(0)=v_0,\quad\dot{v}(0)=v_1,}\ee satisfies
\be\label{3.17n}v(T)=\dot{v}(T)=0,\ee where $\tilde{\X}^m_0(T)$ is
the Banach space of all function with (\ref{3.2n}) and the norm
$(\ref{3.4n})$ but with  (\ref{3.3n}) replaced by
\be\label{3.3!n}\var^{(k)}(T)=0,\quad x\in\Ga_0,\quad0\leq k\leq
m-1.\ee

Then Theorem 1.5 follows by the following

\begin{thm}
Let an equilibrium $w\in H^{m+1}(\Om)\cap H^1_{\Ga_1}(\Om)$ be
exactly controllable. Then there exists a $T_0>0$ such that for any
$T>T_0$ and
$$(v_0,v_1)\in \e(H^{m+1}(\Om)\cap H^1_{\Ga_1}(\Om)\i)\times
\e(H^{m}(\Om)\cap H^1_{\Ga_1}(\Om)\i),$$ there is a
$\var\in\tilde{\X}^m_{0 N}(T)$ such that the solution
$$v\in\cap_{k=0}^mC^k\e([0,T], H^{m-k}(\Om)\i)$$ of the problem
$(\ref{3.15n})$ satisfies $(\ref{3.17n})$.
\end{thm}

As in Section 3, we shall work out the smooth control from the
distributed control theory.

We start with the dual system of the problem (\ref{3.15n})
\be\label{16n} \cases{\ddot{\phi}=\div A(x,\nabla w)\nabla
\phi,\quad (t,x)\in(0,T)\times\Om,\cr
\phi|_{\Ga_1}=\phi_{\nu_A}|_{\Ga_0}=0,\quad t\in(0,T),\cr
\phi(0)=\phi_0,\quad\dot{\phi}(0)=\phi_1.}\ee

We shall need the following observability inequality to get rid of
a lower order term in Lemma 4.4 later: There exists a $T_1>0$ such
that for any $T>T_1$, there is a constant $c_T>0$ for which
\be\label{17n} c_T\int_{\wp_0}\dot{\phi}^2d\wp\geq
E(\phi_0,\phi_1),\ee where $\phi$ is the solution of the problem
(\ref{16n}) and
$$\wp_0=(0,T)\times\Ga_0,\quad
E(\phi_0,\phi_1)=\|A^{1/2}\nabla\phi_0\|^2+\|\phi_1\|^2,\quad
(\phi_0,\phi_1)\in H^1_{\Ga_1}(\Om)\times L^2(\Om),$$ whenever the
left-hand side is finite.

The inequality (\ref{17n}) was established by Lasiecka and Triggiani
\cite{LT} for the classical Laplacian where $\div A(x,\nabla
w)\nabla\phi=\Delta\phi$ and was extended to the case of the
variable coefficients with a first order term in Lasiecka,
Triggiani, and Yao \cite{LTY}, under some geometrical conditions.

The lemma below follows by Lasiecka, Triggiani, and Yao \cite{LTY},
Theorem 3.2, where a uniqueness result, needed, is given by
Triggiani and Yao \cite{TY}, Theorem 10.1.1.

\begin{lem}
Let $w\in H^{m+1}(\Om)\cap H^1_{\Ga_1}(\Om)$ be exactly controllable
such that the assumption $(\ref{2.3*})$ is true and let $\Ga_1$ be
such that $(\ref{2.9*})$ holds. There exists a $T_1>0$ such that for
any $T>T_1$, there is a constant $c_T>0$ for which the inequality
$(\ref{17n})$ is true whenever the left-hand side is finite.
\end{lem}

However, to find  out the smooth control, one-side observability
estimates, as in (\ref{17n}), are insufficient. We have to seek to
establish boundary estimates of another type controlled by the
initial energy both sides from below and also from above, as in
(\ref{3.21}).

Let $\varepsilon>0$ be given small. Let $\eta_\varepsilon\in
C^\infty(\R)$ be such that $0\leq\eta_\varepsilon\leq1$ and
$$\eta_\varepsilon(t)=1 \quad t\leq-\varepsilon;\quad
\eta_\varepsilon(t)=0\quad t\geq0.$$ For any $T>\varepsilon$, let
\be\label{34n}z(t)=\eta_\varepsilon(t-T).\ee Then
$$z(t)=1\quad 0\leq t\leq T-\varepsilon;\quad z(t)=0\quad t\geq T.$$

\begin{lem}
Let \be\label{24n}g=A^{-1}(x,\nabla w)\ee be the Riemannian metric on $\ol{\Om}$.
Let $\phi$ solve the problem \be\label{19n}
\ddot{\phi}=\div A(x,\nabla w)\nabla\phi\quad (t,x)\in\Sigma,\ee
where $\Sigma=(0,T)\times\Om$. Let $H$ be a vector field on
$\ol{\Om}$ and $P\in C^2(\ol{\Om})$ be a function.
Then \beq
\label{20n}&&\int_\wp z
\e[H(\phi)\phi_{\nu_A}+\frac{1}{2}(\dot{\phi}^2-|\nabla_g\phi|_g^2)\<H,\nu\>\i]d\wp\nonumber\\
&&=-\e(\phi_1,H(\phi_0)\i)-\int_{T-\varepsilon}^T\dot{z}\e(\dot{\phi},H(\phi)\i)dt\nonumber\\
&&\quad+\int_\Sigma
z\e[D_gH(\nabla_g\phi,\nabla_g\phi)+\frac{1}{2}(\dot{\phi}^2-|\nabla_g\phi|^2)\div
H\i]d\Sigma,\eeq where $\wp=(0,T)\times\Ga$, and \beq\label{21n}
&&\int_\Sigma zP\e(\dot{\phi}^2-|\nabla_g\phi|_g^2\i)d\Sigma\nonumber\\
&&=-\e(\phi_1,P\phi_0\i)-\int_{T-\varepsilon}^T\dot{z}(\dot{\phi},P\phi)dt-\frac{1}{2}\int_\Sigma
z\phi^2\A Pd\Sigma\nonumber\\
&&\quad+\int_\wp
z\e[\frac{1}{2}\phi^2P_{\nu_A}-P\phi\phi_{\nu_A}\i]d\wp.\eeq
\end{lem}

{\bf Proof.}\,\,\,We multiply the equation (\ref{19n}) by $zH(\phi)$
and $zP\phi$, respectively, integrate by parts over
$\Sigma=(0,T)\times\Om$,
 and obtain the identities (\ref{20n}) and (\ref{21n}),
see Yao \cite{Y}, Proposition 2.1.  {\bf\textbf{$\|$}}\\

Let \be\label{22n}
\Psi(\var,\phi)=\int_{\wp_0}z\e(\dot{\var}\dot{\phi}-\<\nabla_{\Ga_g}\var,\nabla_{\Ga_g}\phi\>_g\i)h_0d\wp,\ee
where $\var$ and $\phi$ solve the problem (\ref{16n}) with the
initial data $(\var_0,\var_1)$ and $(\phi_0,\phi_1)$, respectively,
and \be\label{23n}h_0=\<H_0,\nu\>,\quad
H_0=2\rho_g\nabla_g\rho_g,\quad x\in \Ga,\ee and $\rho_g$ is the
distance function of the Riemannian metric $g$ in (\ref{24n}).

The second observability estimate we need is the following

\begin{lem}
Let $w\in H^{m+1}(\Om)\cap H^1_{\Ga_1}(\Om)$ be exactly controllable
such that the assumption $(\ref{2.3*})$ is true and let $\Ga_1$ be
such that $(\ref{2.9*})$ holds. Let $\varepsilon>0$ be given small.
There are constant $c_{\varepsilon1}>0$, $c_{\varepsilon2}>0$, and
$c_0>0$, independent of time $t$ and solutions $\phi$ of the problem
$(\ref{16n})$, such that for any $T>\varepsilon$
$$c_{\varepsilon2}TE(\phi_0,\phi_1)\geq\Psi(\phi,\phi)+c_0
\int_{\wp_0}\phi^2d\wp+c_0\int_\Sigma\phi^2d\Sigma\geq\e[\rho_0(T-\varepsilon)
-c_{\varepsilon1}\i]E(\phi_0,\phi_1),$$ \be\label{25n}\mbox{for
all}\quad (\phi_0,\phi_1)\in H^1_{\Ga_1}(\Om)\times L^2(\Om),\ee
where $\rho_0>0$ is given in $(\ref{2.3*})$.
\end{lem}

{\bf Proof.}\,\,\,We take $P=\div H_0-\rho_0$ in the identity
(\ref{21n}) and obtain the estimate\beq\label{26n}
&&\int_\Sigma z(\div H_0-\rho_0)\e(\dot{\phi}^2-|\nabla_g\phi|_g^2\i)d\Sigma\nonumber\\
&&\geq-c_{\varepsilon1}E(\phi_0,\phi_1)-c_0\int_\Sigma
\phi^2\Sigma-c_0\int_\wp z\phi^2\wp,\eeq where the boundary
conditions $\phi|_{\Ga_1}=\phi_{\nu_A}|_{\Ga_0}=0$ are used.

Let us take $H=H_0$ in the identity (\ref{20n}) to check the
boundary terms on the left-hand side of the identity (\ref{20n}). On
$\Ga_1$, $\phi_{\Ga_1}=0$ implies $$
H_0(\phi)=\frac{h_0}{|\nu_A|^2_g}\phi_{\nu_A},\quad
|\nabla_g\phi|^2_g=\frac{1}{|\nu_A|^2_g}\phi_{\nu_A}^2,$$ which
implies with $h_0\leq 0$ for $ x\in\Ga_1$ together that
\be\label{27n} \int_{\wp_1}
z\e[H_0(\phi)\phi_{\nu_A}+\frac{1}{2}(\dot{\phi}^2-|\nabla_g\phi|_g^2)h_0\i]d\wp
=\frac{1}{2}\int_{\wp_1}z\frac{\phi_{\nu_A}^2}{|\nu_A|_g^2}h_0d\wp\leq0;\ee
On $\wp_0=(0,T)\times\Ga_0$, $\phi_{\nu_A}=0$ implies $
\nabla_g\phi=\nabla_{\Ga_g}\phi$. We then have via the identity
(\ref{20n}) where $H=H_0$ and (\ref{26n})-(\ref{27n}), (\ref{2.3*}),
that \beq\label{28n} \Psi(\phi,\phi)&&\geq 2\rho_0\int_\Sigma
z|\nabla_g\phi|^2_gd\Sigma-c_{\varepsilon1}E(\phi_0,\phi_1)\nonumber\\
&&\quad+\int_\Sigma z(\dot{\phi}^2-|\nabla_g\phi|^2)\div H_0
d\Sigma\nonumber\\
&&\geq \rho_0\int_\Sigma z(\dot{\phi}^2+|\nabla_g\phi|^2_g)d\Sigma-
c_{\varepsilon1}E(\phi_0,\phi_1)\nonumber\\
&&\quad+\int_\Sigma z(\dot{\phi}^2-|\nabla_g\phi|^2)(\div
H_0-\rho_0)
d\Sigma\nonumber\\
&&\geq[\rho_0(T-\varepsilon)-c_{\varepsilon1}]E(\phi_0,\phi_1)\nonumber\\
&&\quad-c_0\int_\Sigma\phi^2d\Sigma-c_0\int_{\wp_0}\phi^2d\wp_0.\eeq

On the other hand, since $\ol{\Ga}_1\cap\ol{\Ga}_0=\emptyset$, we
take two open sets $\aleph_0$ and $\aleph_1$ in $\R^n$ such that
$\aleph_0\cap\aleph_1=\emptyset$ and $\Ga_i\subset\aleph_i$ for
$i=0$, $1$, respectively. Let $h\in C^\infty(\R^n)$ be such that
$$h(x)=1\quad x\in\aleph_0;\quad h(x)=0\quad x\in\aleph_1.$$
Letting $H=hH_0$ in (\ref{20n}) yields \be\label{29n}
\Psi(\phi,\phi)\leq c_{\varepsilon2}T E(\phi_0,\phi_1).\ee

The lemma follows by the inequalities (\ref{28n}) and (\ref{29n}).
{\bf\textbf{$\|$}}\\

We introduce an operator by \be\label{35n}\A_0 v=\div A(x,\nabla
w)\nabla v,\quad D(\A_0)=\{\,v\in H^2(\Om)\cap H^1_{\Ga_1}(\Om),\,
v_{\nu_A}|_{\Ga_0}=0\,\}.\ee Let $(\phi_0,\phi_1)\in D(\A_0)\times
H^1_{\Ga_1}(\Om)$. Then $(\phi_1,\A_0\phi_0)\in
H^1_{\Ga_1}(\Om)\times L^2(\Om)$. Since $\dot{\phi}$ solves the
problem (\ref{16n}) with the initial data $(\phi_1,\A_0\phi_0)$, the
inequality (\ref{25n}) implies\beq\label{30n}c_{\varepsilon2}T
E(\phi_1,\A_0\phi_0)&&\geq\Psi(\dot{\phi},\dot{\phi})+c_0
\int_{\wp_0}\dot{\phi}^2d\wp+c_0\int_\Sigma\dot{\phi}^2d\Sigma\nonumber\\
&&\geq\e[\rho_0(T-\varepsilon)
-c_{\varepsilon1}\i]E(\phi_1,\A_0\phi_1).\eeq

Let $T_1$  and be given by Lemma 4.1. We fix $T_2>T_1$. Let
$c_{T_2}$ be given by Lemma 4.1. It follows from Lemma 4.1 that for
any $T>T_2+\varepsilon$ \be\label{31n}
\int_\Sigma\dot{\phi}^2d\Sigma\leq T E(\phi_0,\phi_1)\leq T
c_{T_2}\int_0^{T_2}\int_{\Ga_0}\dot{\phi}^2d\wp\leq T
c_{T_2}\int_{\wp_0}z\dot{\phi}^2d\wp.\ee

We introduce a bilinear form by \be\label{32n}
\Psi_*(\var,\phi)=\Psi(\var,\phi)+c_T\int_{\wp_0}z\var\phi d\wp,\ee
where \be\label{75n}c_T=c_0(1+Tc_{T_2}).\ee

Then the inequalities (\ref{30n}) and (\ref{31n}) yield

\begin{lem} For any $T>T_2+\varepsilon$ and $(\phi_0,\phi_1)\in D(\A_0)\times
H^1_{\Ga_1}(\Om)$,\be\label{33n} c_{\varepsilon2}T
E(\phi_1,\A_0\phi_0)\geq\Psi_*(\dot{\phi},\dot{\phi})
\geq\e[\rho_0(T-\varepsilon)
-c_{\varepsilon1}\i]E(\phi_1,\A_0\phi_0).\ee
\end{lem}

We now go back to the control problem in Theorem 4.1.

Given $(\phi_0,\phi_1)\in D(\A_0)\times H^1_{\Ga_1}(\Om)$, the
problem (\ref{16n}) admits a unique solution. We then solve the
problem \be\label{3.19n} \cases{\ddot{\psi}=\div A(x,\nabla w)\nabla
\psi,\quad (t,x)\in(0,T)\times\Om,\cr\psi(T)=\dot{\psi}(T)=0,\quad
x\in\Om, \cr \psi|_{\Ga_1}=0,\cr
\psi_{\nu_A}|_{\Ga_0}=z\e[(\phi^{(3)}-\Delta_{\Ga_g}\dot{\phi})h_0
-\lam_T\dot{\phi}\i],\quad t\in(0,T),}\ee where $\phi$ is produced
by $(\ref{16n})$, $z$ and $h_0$ are given in (\ref{34n}) and
(\ref{23n}), respectively, and \be\label{74n}
\lam_T=c_T+\frac{1}{2}\sup_{x\in\Ga_0}|\Delta_{\Ga_g}h_0|,\ee and
$c_T$ is given by (\ref{75n}).

We define $\Lam_N$: $ D(\A_0)\times H^1_{\Ga_1}(\Om)\rw
\e(H^1_{\Ga_1}(\Om)\i)'\times L^2(\Om)$ by \be\label{3.25n}
\Lam_N(\phi_0, \phi_1)=(\dot{\psi}(0), -\psi(0)),\ee where
$\e(H^1_{\Ga_1}(\Om)\i)'$ is the dual space of $H^1_{\Ga_1}(\Om)$.
Let $\var$ solve the problem (\ref{16n}) with the initial data
$(\var_0,\var_1)$. After we multiply (\ref{3.19n}) by $\var$ and
integrate by parts, we obtain \be\label{3.20n} \<\Lam_N(\phi_0,
\phi_1), (\var_0, \var_1)\>_{L^2(\Om)\times
L^2(\Om)}=-\int_{\wp_0}\psi_{\nu_A}\var d\wp.\ee

Let $k\geq 1$ be an integer. Let $\Xi_{0N}^k(\Om)$ consist of the
functions $u$ in $H^k(\Om)$ with the boundary
conditions\be\label{3.24n}\cases{\A_0^iu|_{\Ga_1}=(\A_0^{i}u)_{\nu_A}|_{\Ga_0}=0,\quad
0\leq i\leq l-1,\quad\mbox{if}\quad k=2l;\cr
\A_0^iu|_{\Ga_1}=(\A_0^{j}u)_{\nu_A}|_{\Ga_0}=0,\quad 0\leq i\leq
l,\,\,0\leq j\leq l-1,\quad\mbox{if}\quad k=2l+1,}\ee and with the
norms of $H^k(\Om)$ where $\A_0$ is given by (\ref{35n}).

The smooth controls with the Neumann action are provided by the
following

\begin{lem}
Let $w\in H^{m+1}(\Om)\cap H^1_{\Ga_1}(\Om)$ be exactly controllable
such that the assumption $(\ref{2.3*})$ is true and let $\Ga_1$ be
such that $(\ref{2.9*})$ holds. Let $k\geq1$ be an integer and
$\Lam_N$ be given by $(\ref{3.25n})$ where $\psi$ is the solution of
the problem $(\ref{3.19n})$. Then there exists a $T_0>0$ such that
for any $T>T_0$, there are $c_1>0$ and $c_2>0$, which depend on $T$,
satisfying $$ c_1\|(\phi_0,\phi_1)\|_{H^{k+1}(\Om)\times
H^k(\Om)}\leq\|\Lam_N(\phi_0,\phi_1)\|_{H^{k-2}(\Om)\times
H^{k-1}(\Om)}\leq c_2\|(\phi_0,\phi_1)\|_{H^{k+1}(\Om)\times
H^k(\Om)},$$ \be\label{3.27n}\any\,\, (\phi_0,\phi_1)\in\Xi_{0
N}^{k+1}(\Om)\times\Xi_{0 N}^k(\Om),\ee where for $k=1$,
$H^{k-2}(\Om)\times H^{k-1}(\Om)=\e(H^1_{\Ga_1}(\Om)\i)'\times
L^2(\Om)$. In particular, $\Lam_N$ are isomorphisms from $\Xi_{0
N}^{2}(\Om)\times\Xi_{0 N}^1(\Om)$ onto
$\e(H^1_{\Ga_1}(\Om)\i)'\times L^2(\Om)$ and from $\Xi_{0
N}^{k+1}(\Om)\times\Xi_{0 N}^k(\Om)$ onto $\e(H^{k-2}(\Om)\cap
H^1_{\Ga_1}(\Om)\i)\times \e( H^{k-1}(\Om)\cap H^1_{\Ga_1}(\Om)\i)$
for $k\geq2$, respectively.
\end{lem}

{\bf Proof.}\,\,\,By induction.

Let $k=1$. Let $(\phi_0,\phi_1)\in\aleph_{0 N}^2(\Om)\times
\aleph_{0 N}^1(\Om)$ be given. For any $(\var_0,\var_1)\in\aleph_{0
N}^2(\Om)\times \aleph_{0 N}^1(\Om)$ given, suppose that $\var$
solves the problem (\ref{16n}) with the initial data
$(\var_0,\var_1)$. Then $\dot{\var}$ solves the problem (\ref{16n})
with the initial $(\var_1,\A_0\var_0)\in H^1_{\Ga_1}(\Om)\times
L^2(\Om)$. By the formula (\ref{3.20n}), we obtain \beq\label{36n}
&&\e(\dot{\psi}(0),\var_1\i)-\e(\psi(0),\A_0\var_0\i)\nonumber\\
&&=\<\Lam_N(\phi_0, \phi_1), (\var_1, \A_0\var_0)\>_{L^2(\Om)\times
L^2(\Om)}=-\int_{\wp_0}\psi_{\nu_A}\dot{\var} d\wp\nonumber\\
&&=\Psi_*(\dot{\phi},\dot{\var})+\int_{\Ga_0}h_0\var_1\Delta_{\Ga_g}\phi_0d\Ga+
\int_{T-\varepsilon}^T\int_{\Ga_0}\dot{z}h_0\dot{\var}\ddot{\phi}
d\Ga dt\nonumber\\
&&\quad-\int_{\wp_0}z\dot{\phi}\nabla_{\Ga_g}h_0(\dot{\var})d\wp+
\frac{1}{2}\sup_{x\in\Ga_0}|\Delta_{\Ga_g}h_0|\int_{\wp_0}z\dot{\phi}\dot{\var}d\wp.\eeq

It follows from (\ref{36n}) and Lemma 4.4 that \beq\label{38n}
&&\e(\dot{\psi}(0),\var_1\i)+\e(\psi(0),\A_0\var_0\i)\nonumber\\
&&\leq
\Psi^{1/2}_*(\dot{\phi},\dot{\phi})\Psi^{1/2}_*(\dot{\var},\dot{\var})
+c\|\Delta_{\Ga_g}\phi_0\|_{H^{-1/2}(\Ga_0)}\|\var_1\|_{H^{1/2}(\Ga_0)}\nonumber\\
&&\quad+c_\varepsilon\int_{T-\varepsilon}^T\|\Delta\phi\|_{H^{-1/2}(\Ga_0)}\|\dot{\var}\|_{H^{1/2}(\Ga_0)}dt
+c\int_0^T\|\dot{\var}\|_{H^{1/2}(\Ga_0)}
\|\dot{\phi}\|_{H^{1/2}(\Ga_0)}dt\nonumber\\
&&\leq
\hat{c}_TE^{1/2}(\phi_1,\A_0\phi_0)E^{1/2}(\var_1,\A_0\var_0),\eeq
for any $T>T_2+\varepsilon$ and all $(\var_0,\var_1)\in \aleph_{0
N}^2(\Om)\times \aleph_{0 N}^1(\Om)$, which gives\be\label{39n}
\|\dot{\psi}(0)\|^2_{\e(H^1_{\Ga_1}(\Om)\i)'}+\|\psi(0)\|^2\leq
\hat{c}_TE(\phi_1,\A_0\phi_0),\quad\any\,\,(\phi_0,\phi_1)\in\aleph_{0
N}^2(\Om)\times \aleph_{0 N}^1(\Om).\ee

Furthermore, letting $(\var_0,\var_1)=(\phi_0,\phi_1)$ in the
identity (\ref{36n}) yields, via Lemma 4.4, \beq\label{37n}
&&\e(\|\dot{\psi}(0)\|_{\e(H^1_{\Ga_1}(\Om)\i)'}^2+\|\psi(0)\|^2\i)^{1/2}E^{1/2}(\var_1,\A_0\var_0)\nonumber\\
&&\geq \e[\rho_0(T-\varepsilon)
-c_{\varepsilon1}\i]E(\phi_1,\A_0\phi_0)-c_\varepsilon
E(\phi_1,\A_0\phi_0)\nonumber\\
&&\geq \e[\rho_0(T-\varepsilon)
-c_{\varepsilon1}\i]E(\phi_1,\A_0\phi_0),\eeq for
$T>T_2+\varepsilon$, where the constant $c_{\varepsilon1}$ may be
different from that in Lemma 4.4 but is independent of time $t$ and
solutions $\phi$.

Combining (\ref{39n}) and (\ref{37n}), we have obtained a $T_0>0$
such that the inequality (\ref{3.27n}) is true for $k=1$.

We assume that the inequality (\ref{3.27n}) is true for some
$k\geq1$. We shall prove it holds true with $k$ replaced by $k+1$.

{\bf Case I}\,\,\,Let $k=2l$ for some $l\geq1$. Firstly, we assume
that \be\label{40n}(\phi_0,\phi_1)\in \aleph_{0
N}^{2k+2}(\Om)\times\aleph_{0 N}^{2k+1}(\Om).\ee

Suppose that $\var$ solves the problem (\ref{16n}) with an initial
data $(\var_0,\var_1)\in\aleph_{0 N}^{2k+2}(\Om)\times\aleph_{0
N}^{2k+1}(\Om)$. Then $\var^{(2i)}$ and $\var^{(2j+1)}$ solve the
problem (\ref{16n}) with the initial data
$(\A^i_0\var_0,\A_0^i\var_1)$ and $(\A_0^j\var_1,\A_0^{j+1}\var_0)$,
respectively, for $0\leq i\leq k+1$ and $0\leq j\leq k-1$.

{\bf Step 1}\,\,\,The following identity is true. \beq\label{49n} &&
-\e(A\nabla(\A_0^{l-1}\dot{\psi}(0)),\nabla(\A_0^l\var_1)\i)-\e(\A_0^l\psi(0),\A_0^{l+1}\var_0\i)\nonumber\\
&&=\Psi_*(\phi^{(k+1)},\var^{(k+1)})+\frac{1}{2}\sup_{x\in\Ga_0}|\Delta_{\Ga_g}h_0|\int_{\wp_0}z\phi^{(k+1)}\var^{(k+1)}d\wp\nonumber\\
&&\quad
-\int_{\wp_0}z\var^{(k+1)}\nabla_{\Ga_g}h_0(\phi^{(k+1)})d\wp+
\sum_{j=1}^{k-1}\int_{T-\varepsilon}^T\int_{\Ga_0}z^{(j)}\phi^{(k+2-j)}\var^{(k+2)}h_0d\Ga
dt\nonumber\\
&&\quad\quad+\sum_{j=1}^k\int_{T-\varepsilon}^T\int_{\Ga_0}z^{(j)}\var^{(k+1)}\e(\Delta_{\Ga_g}\phi^{(k+1-j)}h_0+
\lam_T\phi^{(k+1-j)}\i)d\Ga dt.\eeq

{\em Proof of $(\ref{49n})$}\,\,\, Using $\var^{(2k+1)}$ in place of
$\var$ in the formula (\ref{3.20n}), we obtain
\beq\label{41n}&&\e(\dot{\psi}(0),\A^k_0\var_1\i)-\e(\psi(0),\A_0^{k+1}\var_0\i)=-\int_{\wp_0}\psi_{\nu_A}\var^{(2k+1)}d\wp\nonumber\\
&&=-\int_{\wp_0}h_0z\phi^{(3)}\var^{(2k+1)}d\wp+\int_{\wp_0}zh_0\var^{(2k+1)}\Delta_{\Ga_g}\dot{\phi}d\wp+\lam_T
\int_{\wp_0}z\dot{\phi}\var^{(2k+1)}d\wp\nonumber\\
&&=\Term 1+\Term 2+\lam_T\Term 3, \eeq where $\lam_T$ is given by
(\ref{74n}).

We compute the terms in the right-hand side of (\ref{41n}) by
integrating by parts over $\wp_0=(0,T)\times\Ga_0$, respectively, as
\beq\label{42n} \Term
1&&=\int_{\wp_0}z\phi^{(k+2)}\var^{(k+2)}h_0d\wp
+\sum_{j=0}^{k-2}(-1)^j\int_{\Ga_0}\phi^{(j+3)}(0)\var^{(2k-j)}(0)h_0d\Ga
\nonumber\\
&&+
\sum_{j=1}^{k-1}\int_{T-\varepsilon}^T\int_{\Ga_0}z^{(j)}\phi^{(k+2-j)}\var^{(k+2)}h_0d\Ga dt\nonumber\\
&&=\int_{\wp_0}z\phi^{(k+2)}\var^{(k+2)}h_0d\wp+\sum_{j=1}^l\int_{\Ga_0}\A_0^j\phi_1\A_0^{k+1-j}\var_0
h_0d\Ga\nonumber\\
&&\quad-\sum_{j=2}^l\int_{\Ga_0}\A_0^j\phi_0\A_0^{k+1-j}\var_1
h_0d\Ga\nonumber\\
&&+
\sum_{j=1}^{k-1}\int_{T-\varepsilon}^T\int_{\Ga_0}z^{(j)}\phi^{(k+2-j)}\var^{(k+2)}h_0d\Ga
dt; \eeq \beq\label{43n}\Term
2&&=\int_{\wp_0}zh_0\var^{(k+1)}\Delta_{\Ga_g}\phi^{(k+1)}d\wp
+\sum_{j=1}^k\int_{T-\varepsilon}^T\int_{\Ga_0}z^{(j)}\var^{(k+1)}\Delta_{\Ga_g}\phi^{(k+1-j)}h_0d\Ga
dt\nonumber\\
&&\quad+\sum_{j=0}^{k-1}(-1)^{j+1}\int_{\Ga_0}h_0\var^{(2k-j)}(0)\Delta_{\Ga_g}\phi^{(j+1)}(0)d\Ga\nonumber\\
&&=-\int_{\wp_0}z\<\nabla_{\Ga_g}\phi^{(k+1)},\nabla_{\Ga_g}\var^{(k+1)}\>_gh_0d\wp-\int_{\wp_0}z\var^{(k+1)}\nabla_{\Ga_g}h_0(\phi^{k+1)})d\wp\nonumber\\
&&\quad+\sum_{j=1}^l\int_{\Ga_0}h_0\A_0^{k-j}\var_1\Delta_{\Ga_g}\A_0^j\phi_0d\Ga
-\sum_{j=0}^{l-1}\int_{\Ga_0}h_0\A_0^{k-j}\var_0\Delta_{\Ga_g}\A_0^j\phi_1d\Ga\nonumber\\
&&\quad
+\sum_{j=1}^k\int_{T-\varepsilon}^T\int_{\Ga_0}z^{(j)}\var^{(k+1)}\Delta_{\Ga_g}\phi^{(k+1-j)}h_0d\Ga
dt;\eeq \beq\label{44n} \Term
3&&=\int_{\wp_0}z\phi^{(k+1)}\var^{(k+1)}d\wp+\sum_{j=1}^k\int_{T-\varepsilon}^T\int_{\Ga_0}z^{(j)}\phi^{(k+1-j)}\var^{(k+1)}d\Ga
dt\nonumber\\
&&\quad+\sum_{j=1}^l\int_{\Ga_0}\A_0^j\phi_0\A_0^{k-j}\var_1d\Ga-\sum_{j=0}^{l-1}\int_{\Ga_0}\A_0^j\phi_1\A_0^{k-j}\var_0d\Ga.
\eeq

Moreover, via the problem (\ref{3.19n}), we have, on $\Ga_0$ for
$j\geq0$,
\beq\label{47n}&&\e(\A_0^j\dot{\psi}(0)\i)_{\nu_A}=\psi^{(2j+1)}_{\nu_A}(0)
=(\A_0^{j+2}\phi_0-\Delta_{\Ga_g}\A^{j+1}_0\phi_0)h_0-\lam_T\A_0^{j+1}\phi_0;
\eeq
\be\label{48n}\e(\A_0^j\psi(0)\i)_{\nu_A}=(\A_0^{j+1}\phi_1-\Delta_{\Ga_g}\A_0^j\phi_1)h_0-\lam_T\A_0^j\phi_1.\ee

We substitute (\ref{42n})-(\ref{48n}) into (\ref{41n}) to yield
\beq\label{50n}&&\e(\dot{\psi}(0),\A^k_0\var_1\i)-\e(\psi(0),\A_0^{k+1}\var_0\i)=\Psi_*(\phi^{(k+1)},\var^{(k+1)})
\nonumber\\
&&\quad+\frac{1}{2}\sup_{x\Ga_0}|\Delta_{\Ga_g}h_0|\int_{\wp_0}z\phi^{(k+1)}\var^{(k+1)}d\wp\nonumber\\
&&\quad
-\int_{\wp_0}z\var^{(k+1)}\nabla_{\Ga_g}h_0(\phi^{k+1)})d\wp+
\sum_{j=1}^{k-1}\int_{T-\varepsilon}^T\int_{\Ga_0}z^{(j)}\phi^{(k+2-j)}\var^{(k+2)}h_0d\Ga
dt\nonumber\\
&&\quad\quad+\sum_{j=1}^k\int_{T-\varepsilon}^T\int_{\Ga_0}z^{(j)}\var^{(k+1)}\e(\Delta_{\Ga_g}\phi^{(k+1-j)}h_0+\phi^{(k+1-j)}\i)d\Ga
dt\nonumber\\
&&\quad+\sum_{j=0}^{l-1}\int_{\Ga_0}\e(\A_0^j\psi(0)\i)_{\nu_A}\A_0^{k-j}\var_0d\Ga-\sum_{j=0}^{l-2}\int_{\Ga_0}
\e(\A_0^j\dot{\psi}(0)\i)_{\nu_A}\A_0^{k-1-j}\var_1d\Ga.\eeq

On the other hand, for $(\var_0,\var_1)\in\aleph_{0
N}^{2k+2}(\Om)\times\aleph_{0 N}^{2k+1}(\Om)$, we obtain via the
Green
formula\beq\label{45n}\e(\dot{\psi}(0),\A^k_0\var_1\i)&&=-\e(A\nabla(\A_0^{l-1}\dot{\psi}(0)),\nabla(\A_0^l\var_1)\i)\nonumber\\
&&\quad-\sum_{j=0}^{l-2}\int_{\Ga_0}
\e(\A_0^j\dot{\psi}(0)\i)_{\nu_A}\A_0^{k-1-j}\var_1d\Ga; \eeq
\beq\label{46n}
-\e(\psi(0),\A_0^{k+1}\var_0\i)&&=-\e(\A_0^l\psi(0),\A_0^{l+1}\var_0\i)+\sum_{j=0}^{l-1}\int_{\Ga_0}\e(\A_0^j\psi(0)\i)_{\nu_A}\A_0^{k-j}\var_0d\Ga.\eeq

After substituting (\ref{45n}) and (\ref{46n}) into the left-hand
side of the identity (\ref{50n}) and eliminating the same terms from
the both sides, we obtain the identity (\ref{49n}).

{\bf Step 2}\,\,\,We have \be\label{57nn}
\|(\dot{\psi}(0),\psi(0))\|^2_{H^{k-1}(\Om)\times H^{k}(\Om)}\geq
c_1\|(\phi_0,\phi_1)\|^2_{H^{k+2}(\Om)\times H^{k+1}(\Om)},\ee for
all $(\phi_0,\phi_1)\in\aleph_{0 N}^{k+2}(\Om)\times\aleph_{0
N}^{k+1}(\Om)$ and $T$ large.

{\em Proof of $(\ref{57nn})$}\,\,\, Replace $\phi$ with
$\phi^{(k)}$ in the inequality (\ref{33n}) and obtain
\be\label{50n*} c_{\varepsilon2}T
E(\A_0^l\phi_1,\A_0^{l+1}\phi_0)\geq\Psi_*(\phi^{(k+1)},\phi^{(k+1)})
\geq\e[\rho_0(T-\varepsilon)
-c_{\varepsilon1}\i]E(\A_0^l\phi_1,\A_0^{l+1}\phi_0),\ee for all
$(\phi_0,\phi_1)\in\aleph_{0 N}^{k+2}(\Om)\times\aleph_{0
N}^{k+1}(\Om)$.

We let $(\var_0,\var_1)=(\phi_0,\phi_1)$ in the identity (\ref{49n})
and observe that
\be\label{53n}\int_{\wp_0}z\phi^{(k+1)}\nabla_{\Ga_g}(\phi^{(k+1)})d\wp=
-\frac{1}{2}\int_{\wp_0}[\phi^{(k+1)}]^2\Delta_{\Ga_g}h_0d\wp;\ee
\beq\label{52n}
&&\Bigm|\sum_{j=1}^{k-1}\int_{T-\varepsilon}^T\int_{\Ga_0}z^{(j)}\phi^{(k+2-j)}\phi^{(k+2)}h_0d\Ga
dt\Bigm|\nonumber\\
&&\leq c_\varepsilon\sum_{j=1}^{k-2}\sup_{T-\varepsilon\leq t\leq
T}\|\phi^{(k+2-j)}\|_{L^2(\Ga_0)}\|\phi^{(k+1)}\|_{L^2(\Ga_0)}\leq
c_\varepsilon E(\A_0^l\phi_1,\A_0^{l+1}\phi_0);\eeq
\beq\label{54n}&&\Bigm|\sum_{j=1}^k\int_{T-\varepsilon}^T\int_{\Ga_0}z^{(j)}\phi^{(k+1)}\e(\Delta_{\Ga_g}\phi^{(k+1-j)}h_0+
\lam_T\phi^{(k+1-j)}\i)d\Ga dt\Bigm|\nonumber\\
&&\leq c_\varepsilon\sum_{j=1}^{k}\sup_{T-\varepsilon\leq t\leq
T}\|\Delta_{\Ga_g}\phi^{(k+1-j)}\|_{H^{-1/2}(\Ga_0)}\|\phi^{(k+1)}\|_{H^{1/2}(\Ga_0)}\nonumber\\
&&\quad+\lam_Tc_\varepsilon\sum_{j=1}^{k}\sup_{T-\varepsilon\leq
t\leq
T}\|\phi^{(k+1-j)}\|_{H^{1/2}(\Ga_0)}\|\phi^{(k+1)}\|_{H^{-1/2}(\Ga_0)}\nonumber\\
&&\leq c_\varepsilon
E(\A_0^l\phi_1,\A_0^{l+1}\phi_0)+\lam_Tc_\varepsilon
E(\A_0^l\phi_0,\A_0^l\phi_1).\eeq

We then obtain by setting $(\var_0,\var_1)=(\phi_0,\phi_1)$ in
(\ref{49n}) and via (\ref{50n})-(\ref{54n})
\beq\label{51n}&&\|A\nabla(\A_0^{l-1}\dot{\psi}(0))\|^2+
\|\A_0^l\psi(0)\|^2+\lam_Tc_\varepsilon E(\A_0^l\phi_0,\A_0^l\phi_1)\nonumber\\
&&\geq\e[\rho_0(T-\varepsilon)
-c_{\varepsilon1}\i]E(\A_0^l\phi_1,\A_0^{l+1}\phi_0).\eeq Next, the
inductive assumption that the inequality (\ref{3.27}) holds for $k$
implies that for $T$ large there is $c>0$ such that \be\label{56n}
c\|(\dot{\psi}(0),\psi(0))\|^2_{H^{k-2}(\Om)\times H^{k-1}(\Om)}\geq
E(\A_0^l\phi_0,\A_0^l\phi_1).\ee

Combining (\ref{51n}) and (\ref{56n}) yields that the inequality
(\ref{57nn}) is true for all $(\phi_0,\phi_1)\in\aleph_{0
N}^{2k+2}(\Om)\times\aleph_{0 N}^{2k+1}(\Om)$ and $T$ large.
 Since $\aleph_{0
N}^{2k+2}(\Om)\times\aleph_{0 N}^{2k+1}(\Om)$  is dense in
$\aleph_{0 N}^{k+2}(\Om)\times\aleph_{0 N}^{k+1}(\Om)$, then
inequality (\ref{57nn}) is actually true for all
$(\phi_0,\phi_1)\in\aleph_{0 N}^{k+2}(\Om)\times\aleph_{0
N}^{k+1}(\Om)$.

{\bf Step 3}\,\,\,There is $c_2>0$ such that \be\label{62n}
\|\dot{\psi}(0)\|^2_{k-1}+\|\psi(0)\|^2_k\leq
c_2\|(\phi_0,\phi_1)\|^2_{H^{k+2}(\Om)\times H^{k+1}(\Om)},\ee for
all $(\phi_0,\phi_1)\in\aleph_{0 N}^{k+2}(\Om)\times\aleph_{0
N}^{k+1}(\Om)$.

{\em Proof of $(\ref{62n})$}\,\,\,We let $\var_1=0$ and
$\var_0\in\aleph_{0 N}^{2k+2}(\Om)$ in the identity (\ref{49n}) and
use the inequality (\ref{50n}). We obtain\beq\label{63n}
|(\A_0^l\psi(0),\A^{l+1}_0\var_0)|&&\leq\Psi^{1/2}_*(\phi^{(k+1)},\phi^{(k+1)})\Psi^{1/2}_*(\var^{(k+1)},\var^{(k+1)})
\nonumber\\
&&\quad+c\int_0^T\|\phi^{(k+1)}\|_{H^{1/2}(\Ga_0)}\|\var^{(k+1)}\|_{H^{1/2}(\Ga_0)}dt\nonumber\\
&&\quad+c\sum_{j=1}^{k-1}\int_{T-\varepsilon}^T\|\phi^{(k+2-j)}\|_{H^{1/2}(\Ga_0)}\|\var^{(k+2)}\|_{H^{-1/2}(\Ga_0)}dt
\nonumber\\
&&\quad
+c\sum_{j=1}^k\int_{T-\varepsilon}^T\|\var^{(k+1)}\|_{H^{1/2}(\Ga_0)}\|\Delta_{\Ga_g}\phi^{(k+1-j)}\|_{H^{-1/2}(\Ga_0)}dt
\nonumber\\
&&\leq
cE^{1/2}(\A^l_0\phi_1,\A_0^{l+1}\phi_0)E^{1/2}(\A^l_0\var_1,\A_0^{l+1}\var_0)\nonumber\\
&&=c E(\A^l_0\phi_1,\A_0^{l+1}\phi_0)\|\A^{l+1}_0\var_0\|, \eeq
since $\var_1=0$. Because $\A_0^{l+1}$: $\aleph_{0 N}^{k+2}(\Om)\rw
L^2(\Om)$ is an isomorphism, it follows from (\ref{63n})
that\be\label{64n}\|\A_0^l\psi(0)\|^2\leq c
E((\A^l_0\phi_1,\A_0^{l+1}\phi_0),\ee for all
$(\phi_0,\phi_1)\in\aleph_{0 N}^{k+2}(\Om)\times\aleph_{0
N}^{k+1}(\Om)$.

Next, by the ellipticity of the operator $\A_0$ and the equation in
(\ref{16n}), we have\beq\label{58n}&&\|\psi(0)\|^2_k\leq
c\|\A_0\psi(0)\|^2_{k-2}+c\|\psi_{\nu_A}(0)\|_{k-3/2,\Ga_0}^2+c\|\psi(0)\|^2_{k-1}\nonumber\\
&&\leq
c\|\A_0^2\psi(0)\|^2_{k-4}+c\|\ddot{\psi}_{\nu_A}(0)\|^2_{k-7/2,\Ga_0}+c\|\psi_{\nu_A}(0)\|_{k-3/2,\Ga_0}^2
+c\|\psi(0)\|^2_{k-1}.\nonumber\eeq Repeating this process
gives\be\label{59n}\|\psi(0)\|^2_k\leq
c\|\A_0^l\psi(0)\|^2+c\sum_{j=0}^{l-1}\|\psi_{\nu_A}^{(2j)}(0)\|^2_{k-2j-3/2,\Ga_0}+c\|\psi(0)\|^2_{k-1}.\ee

We use the boundary control of (\ref{3.19n}) and the equation in the
problem (\ref{16n}). We obtain\beq\label{60n}
\sum_{j=0}^{l-1}\|\psi_{\nu_A}^{(2j)}(0)\|^2_{k-2j-3/2,\Ga_0}&&\leq
c\sum_{j=0}^{l-1}\e(\|\phi^{(2j+3)}(0)\|^2_{k-2j-3/2,\,\Ga_0}+
\|\Delta_{\Ga_g}\phi^{(2j+1)}(0)\|^2_{k-2j-3/2,\Ga_0}\i)\nonumber\\
&&\quad+c\sum_{j=0}^{l-1}\|\phi^{(2j+1)}(0)\|^2_{k-2j-3/2,\,\Ga_0}\nonumber\\
&&\leq c\sum_{j=0}^{l-1}\e(\|\phi^{(2j+3)}(0)\|^2_{k-2j-1}+
\|\A_0\phi^{(2j+1)}(0)\|^2_{k-2j-1}\i)\nonumber\\
&&\leq c\|\phi^{(2l+1)}(0)\|^2_1\leq c
E(\A_0^l\phi_1,\A_0^{l+1}\phi_0).\eeq

We combine (\ref{64n})-(\ref{60n}) and use the inductive assumption
$$\|\psi(0)\|^2_{k-1}\leq c\|(\phi_0,\phi_1)\|^2_{H^{k+1}(\Om)\times
H^k(\Om)},$$ and we have \be\label{65n}\|\psi(0)\|^2_k\leq c
E(\A^l_0\phi_1,\A_0^{l+1}\phi_0)\leq
c\|(\phi_0,\phi_1)\|^2_{H^{k+2}(\Om)\times H^{k+1}(\Om)},\ee for all
$(\phi_0,\phi_1)\in\aleph_{0 N}^{k+2}(\Om)\times\aleph_{0
N}^{k+1}(\Om)$.

A similar argument establishes the estimate for $\dot{\psi}(0)$.

{\bf Case II}\,\,\,Let $k=2l+1$ for some $l\geq0$. A similar
argument shows that the inequality (\ref{3.27n}) holds with $k$
replaced by $k+1$ if it is true for $k$.

Finally, the lemma follows by induction.
 {\bf\textbf{$\|$}}\\

We consider the regularity of the control function in the problem
(\ref{3.19n}). Since $\dot{\phi}$ is a lower order term in the the
boundary control of (\ref{3.19n}),
$\phi^{(3)}-\Delta_{\Ga_g}\dot{\phi}$ is the principle part of the
control. The following lemma relates to the regularity of this
principle part.

\begin{lem}
Let $\phi$ solve the problem $(\ref{16n})$ with the initial data
$(\phi_0,\phi_1)\in D(\A_0)\times H^1_{\Ga_1}(\Om)$.
Then\be\label{66n} \ddot{\phi}-\Delta_{\Ga_g}\phi\in L^2(\Ga_0).\ee
Furthermore, if $(\phi_0,\phi_1)\in \aleph_{0 N}^3(\Om)\times
\aleph_{0 N}^2(\Om)$, then \be\label{72n}
\ddot{\phi}-\Delta_{\Ga_g}\phi\in C\e([0,T],H^{1/2}(\Ga_0)\i)\cap
H^1\e((0,T)\times\Ga_0\i).\ee
\end{lem}

{\bf Proof.}\,\,\,Let the Riemann metric $g$ be given by
(\ref{24n}). Then \be\label{67n} \A_0\phi=\Delta_g\phi+F(\phi)\quad
x\in\Om,\ee where $\Delta_g$ is the Laplacian of the metric $g$ and
$F$ is a vector field on $\Om$ give by $$F=\frac{1}{2G}A(x,\nabla
w)\nabla G,$$ $G$ being the determinant of $A^{-1}(x,\nabla w)$.

Using the boundary condition $\phi_{\nu_A}=0$ on $\Ga_0$ and the
relation (\ref{67n}), we obtain\be\label{68n}
\ddot{\phi}-\Delta_{\Ga_g}\phi=\frac{1}{|\nu_A|^2_g}D^2_g\phi(\nu_A,\nu_A)+\<F,\nabla_{\Ga_g}\phi\>_g,\quad
x\in\Ga_0,\ee where $D^2_g\phi(\cdot,\cdot)$ is the Hessian of
$\phi$ in the metric $g$. Since
$\|\nabla_{\Ga_g}\phi\|^2_{L^2(\Ga_0)}\leq c E(\phi_1,\A_0\phi_0)$,
to get the relation (\ref{66n}) it will suffice to
prove\be\label{69n} D^2_g\phi(\nu_A,\nu_A)\in L^2(\Ga_0).\ee

Let $H$ be a vector field on $\ol{\Om}$ such that $$H=0,\quad
x\in\Ga_1;\quad H=\nu_A,\quad x\in\Ga_0.$$ We set
\be\label{71n}\var=H(\phi),\quad x\in\Om.\ee It is easy to check
that $\var$, given by (\ref{71n}), solves the problem with the
Dirichlet boundary
conditions\be\label{70n}\cases{\ddot{\var}=\A_0\var+[H,\A_0]\phi,\quad
(t,x)\in(0,T)\times\Om,\cr\var|_{\Ga}=0,\quad
t\in(0,T),\cr\var(0)=H(\phi_0),\quad\dot{\var}(0)=H(\phi_1).}\ee In
addition, $(\phi_0,\phi_1)\in D(\A_0)\times H^1_{\Ga_1}(\Om)$
implies $(\var(0),\dot{\var}(0))\in H^1_0(\Om)\times L^2(\Om)$. We
use lemma 3.1 to obtain
$$\var_{\nu_A}=D_g^2\phi(\nu_A,\nu_A)+\<\nabla_{\Ga_g}\phi,\,\,(D_g)_{\nu_A}\nu_A\>_g\in
L^2(\Ga_0),$$ which gives the relation (\ref{69n}).

Next, we assume that $(\phi_0,\phi_1)\in \aleph_{0 N}^3(\Om)\times
\aleph_{0 N}^2(\Om)$. Then $(\var(0),\dot{\var}(0))\in
\e(H^2(\Om)\times H^1_0(\Om)\i)\times H^1_0(\Om)$, where $\var$ is
given by (\ref{71n}). A similar argument as in the proof of Lemma
3.3 shows that $$\var_{\nu_A}\in C\e([0,T],H^{1/2}(\Ga)\i)\cap
H^1\e((0,T)\times\Ga\i),$$ which implies that the relation (\ref{72n}) is true.  {\bf\textbf{$\|$}}\\

If  $(\phi_0,\phi_1)\in \aleph_{0 N}^3(\Om)\times \aleph_{0
N}^2(\Om)$, the relation (\ref{72n}) shows that we can find a
control $\var$ in $L^2\e((0,T)\times\Ga_0\i)$ to move one state to
another in the space $L^2(\Om)\times H^1_{\Ga_0}(\Om)$ by the
control scheme in (\ref{3.19n}).\\

{\bf The Proof of Theorem 4.1}\,\,\,By Lemma 4.5, there is $T_0>0$
such that for any $T>T_0$, $\Lam_N$: $\aleph_{0
N}^{m+3}(\Om)\times\aleph_{0 N}^{m+2}(\Om)\rw H^m(\Om)\times
H^{m+1}(\Om)$ is an isomorphism. For any $(v_0,v_1)\in
H^{m+1}(\Om)\times H^m(\Om)$ there is a unique $(\phi_0,\phi_1)\in
\aleph_{0 N}^{m+3}(\Om)\times\aleph_{0 N}^{m+2}(\Om)$ such that the
solution of the problem (\ref{3.15n}) satisfies (\ref{3.17n}) under
the control action
\be\label{73n}\var=z\e[(\phi^{(3)}-\Delta_{\Ga_g}\dot{\phi})h_0
-\lam_T\dot{\phi}\i],\quad x\in\Ga_0,\ee where $\phi$ solves the
problem (\ref{16n}).

To complete the proof, we need to verify $\var\in \tilde{\X}_{0
N}^m(T)$. Indeed, the relation  $(\phi_0,\phi_1)\in \aleph_{0
N}^{m+3}(\Om)\times\aleph_{0 N}^{m+2}(\Om)$ implies that
$(\phi^{(m)}(0), \phi^{(m+1)}(0))\in \aleph_{0
N}^{3}(\Om)\times\aleph_{0 N}^{2}(\Om)$. It follows from Lemma 4.6
that $\phi^{(m+2)}-\Delta_{\Ga_g}\phi^{(m)}\in
C\e([0,T],H^{1/2}(\Ga)\i)\cap H^1\e((0,T)\times\Ga\i)$ which yields
 $\var\in \tilde{\X}_{0
N}^m(T)$.

\def\theequation{5.\arabic{equation}}
\setcounter{equation}{0}
\section{Globally exact controllability; Geometrical conditions }
\hskip\parindent {\bf The Proof of Theorem 1.3}\,\,\,By Theorem 1.2
and the compactness principle it will suffice to prove that $w_\a$:
$[0,1]\rw H^m(\Om)$ is continuous in $\a\in[0,1]$.

It is readily seen that $v_\a=\frac{\pl}{\pl\a}w_a$ is the solution
of the following linear, elliptic problem\be\label{4.1}
\cases{\sum_{ij=1}^na_{ij}(x,\nabla w_\a)v_{\a
x_ix_j}+\sum_{l=1}^n\e[\sum_{ij=1}^na_{ijy_l}(x,\nabla w_\a)w_{\a
x_ix_j}+b_{y_l}(x,\nabla w_\a)\i]v_{\a x_l}=0, \cr
v_\a|_\Ga=w|_{\Ga},}\ee for each $\a\in[0,1]$, and, in addition, by
the maximum principle for the above problem
(\ref{4.1}),\be\label{4.4}\sup_{x\in\Om}|\frac{\pl}{\pl\a}w_{\a}|\leq
\sup_{x\in\Ga}|w|.\ee

Let \be\label{4.2} B(\a)v=\sum_{ij=1}^na_{ij}(x,\nabla
w_\a)v_{x_ix_j},\quad v\in H^2(\Om),\quad \a\in[0,1].\ee
 By the uniform bound (\ref{2.31*}), the ellipticity of the operator $B(\a_0)$, and the estimate (\ref{4.4}),
 we have\beq\label{4.3}&&\|w_\a-w_{\a_0}\|_m\nonumber\\
 &&\leq
c\|B(\a_0)(w_\a-w_{\a_0})\|_{m-2}+c|\a-\a_0|\|w\|_{m-1/2,\,\Ga}+c\|w_\a-w_{\a_0}\|\nonumber\\
&&\leq
c\|B(\a_0)(w_\a-w_{\a_0})\|_{m-2}+c|\a-\a_0|\e(\|w\|_{m-1/2,\,\Ga}+\sup_{x\in\Ga}|w|\i).\eeq

Next, let us estimate $\|\B(\a_0)(w_\a-w_{\a_0})\|_{m-2}$.

$[B(\a_0)-B(\a)]w_\a$ and $b(x,\nabla w_\a)-b(x,\nabla w_{\a_0})$
can be written as  sums of some terms of the form, respectively,
$$f(x,\nabla w_\a,\nabla w_{\a_0})(w_{\a_0x_l}-w_{\a x_l})w_{\a x_ix_j}.$$ Applying
the estimate (\ref{2.7}) to the above products gives, via the bound
(\ref{2.31*}) and the estimate (\ref{4.4}),
\beq\label{4.5}\|B(\a_0)(w_\a-w_{\a_0})\|_{m-2}&&\leq
\|\e(B(\a_0)-B(\a)\i)w_{\a}\|_{m-2}+\|b(x,\nabla w_a)-b(x,\nabla
w_{\a_0})\|_{m-2}\nonumber\\
&&\leq c\|w_\a-w_{\a_0}\|_{m-1}\nonumber\\
&&\leq\varepsilon\|w_\a-w_{\a_0}\|_{m-2}+c_\varepsilon|\a-\a_0|\sup_{x\in\Ga}|w|.\eeq

We obtain the desired result after substituting (\ref{4.5}) into
(\ref{4.3}). {\bf\textbf{$\|$}}\\

{\bf The Proof of Theorem 1.6}\,\,\,The same argument as above
completes the proof.  {\bf\textbf{$\|$}}\\

To end this paper, we prove Proposition 1.1.

{\bf The Proof of Proposition 1.1}\,\,\,We only need to prove the
case of $\kappa>0$. By Yao \cite{Y}, Corollary 1.2, if there are
$x_0\in\ol{\Om}$ and $\gamma>0$ such that \be\label{4.6} \Om\subset
B_{g_w}(x_0,\gamma), \quad 4\gamma^2\kappa<\pi^2,\ee where
$$B_{g_w}(x_0,\gamma)=\{\,x\,|\,x\in\R^n,\,\rho_{g_w}(x_0,x)<\gamma\,\},$$
then the inequality (\ref{2.3*}) is true. To complete the proof, it
will suffice to prove that the condition (\ref{2.34*}) implies
(\ref{4.6}). By (\ref{2.34*}), there is a
$0<\gamma_1<\lam\pi/(2\sqrt{\kappa})$ such that\be\label{4.8}
\Om\subset B(x_0,\gamma_1).\ee For $x\in B(x_0,\gamma_1)$ be given,
$r(t)=tx_0+(1-t)x$ is a curve in $(\R^n,g_w)$ for $0\leq t\leq1$
which connects the points $x_0$ and $x$. Then \beq\label{4.9}
\rho_{g_w}(x_0,x)&&\leq
\int_0^1|\dot{r}(t)|_{g_w}dt=\int_0^1\<A^{-1}(x,\nabla
w)\dot{r}(t),\,\,\dot{r}(t)\>^{1/2}dt\nonumber\\
&&\leq\frac{1}{\lam}|x-y|\leq \frac{\gamma_1}{\lam},\eeq which
implies that (\ref{4.6}) is true with $\gamma=\gamma_1/\lam$.
 {\bf\textbf{$\|$}}\\

{\bf Acknowledgments}\\

The work was in part done while the author was visiting Imperial
College. Hospitality and support by the NSFC China-the Royal
Society joint project is acknowledged and greatly appreciated.
Particular thanks are extended to Prof. George Weiss for many
stimulating discussions. The content of this paper was presented
at a seminar at the University of Donbei.


\begin{thebibliography}{}

\bibitem{B} M. Berger, Nonlinearity and Functional Analysis,
Academic Press, 1977.

\bibitem{BLR} C. Bardos, G. Lebeau, and J. Rauch, Sharp sufficient
conditions for the observation, control and stabilization of waves
from the boundary, {\em SIAM J. Control Optim.}, {\bf 30} (1992),
1024-1065.

\bibitem{CZ} C. Castro, and E. Zuazua, Concentration and lack
observability od waves in highly heterogeneous media, {\em Arch.
Ration. Anal.}, {\bf 164} (2002), no. {\bf 1}, 39-72.

\bibitem{C} M. Cirina, Boundary controllability of nonlinear hyperbolic systems, SIAM
{\em J. Control}, 7(1969), 198-212.

\bibitem{C1} S. Chai, Y. Guo, and P.F. Yao, Boundary feedback
stablization of shallow shells, {\em SIAM J. Control Optim.} {\bf
42}(2003), no. 1, 239-259.

\bibitem{C2} S. Chai and P.F. Yao, Observability inequalities for
thin shells, {\em Science in China (Series A)}, Vol. {\bf 46}, No.
3, 300-311.

\bibitem{DH} C. M. Dafermos and W. J. Hrusa, Energy methods for
quasilinear hyperbolic initial-boundary value problems.
Applications to elastodynamics, {\em Arch. Ration. Mech. Anal.}
{\bf 87} (1985), 267-292.

\bibitem{E} Yu. V. Egorov, Some problems in the theory of optimal
control, {\em Z. Vycisl Mat. i Mat. Fiz.} (1963), no. {\bf 5},
887-904.

\bibitem{GT} D. Gilbarg and N.S. Trudinger, Elliptic Partial
Differential Equations of Second Order, Second Edition and Revised
Third Printing, Springer-Verlag, 1998.

\bibitem{GU} R. Gulliver, I. Lasiecka, W. Littman, and R. Triggiani,
The case for differential geometry in the control of single and
coupled PDEs: the structural acoustic chamber. Geometric methods in
inverse problems and PDE control, 73--181, {\em IMA Vol. Math.
Appl.}, {\bf 137}, Springer, New York, 2004.

\bibitem{Fa} H. O. Fattorini, Boundary control of temperature
distributions in parallepipedon, {\em SIAM J. Control}, {\bf 13}
(1975), no. {\bf 1}, 1-13.

\bibitem{H} L.F. Ho, Observabilit$\acute{e}$ fronti$\acute{e}$re de
l'e'quation des ondes, {\em C.R. Acad. Sci.  Paris S$\acute{e}$r. I
Math.}, {\bf 302}(1986), pp 443-446.

\bibitem{LT} I. Lasiecka and R. Triggiani, Exact controllability of
the wave equation with Neumann boundary control, {\em Appl. Math.
Optimiz.}, {\bf 19} (1989), 243-209.

\bibitem{LT1} I. Lasiecka and R. Triggiani, Uniform
stabilization of a shallow shell model with nonlinear boundary
feedbacks, {\em J. Math. Anal. Appl.}, {\bf  269} (2002), no. 2,
642--688.

\bibitem{LTY} I. Lasiecka, R. Triggiani, and P.F. Yao,
Inverse/observability estimates for second-order hyperbolic
equations with variable systems, {\em J. Math. Anal. Appl.}, {\bf
235} (1999) 13-57.

\bibitem{LR} T.T. Li and B.P. Rao, Exact boundary controllability
for quasilinear hyperbolic systems, {\em SIAM J. Control Optim.}
{\bf 41}(2003), no. 6, 1748-1755.

\bibitem{L} J.L. Lions, Exact controllability, stabilization and
perturbations for distributed system, {\em SIAM Reviews}, Vol {\bf
30}(1988)1-68.

\bibitem{Ru} D. L. Russell, Controllability and stability theory
for linear partial differential equations, Reccent progress and
open questions, {\em SIAM Review}, {\bf 20} (1978), no. {\bf 4},
639-739.

\bibitem{S} E.J.P.G. Schmidt, On a non-linear wave equation and
th e control of an elastic string from one equilibrium location to
another, {\em J. Math. Anal. Appl.} {\bf 272}(2002)536-554.

\bibitem{Se} T. I. Seidman, Two results on exact boundary
controllability of parabolic equations, {\em Applied Math. and
Optimization}, {\bf 11} (1984), no. {\bf 2}, 145-152.

\bibitem{Ta} D. Tataru, Boundary controllability for conservative
PDEs, {\em Appl. Math. Optim.}, {\bf 31} (1995), 257-295.

\bibitem{T} M. E. Taylor, Partial Differential Equations I, Springer-Verlag, 1996.

\bibitem{TY} R. Triggiani and P.F. Yao, Carleman estimate with no
lower-order terms for general Riemann wave equation. Global
uniqueness and observability in one shot, {\em Appl. Math. Optim.},
{\bf 46} (2002) 331-375.

\bibitem{Y}P.F. Yao, On the observability inequalities for the
exact controllability of the wave equation with variable
coefficients, {\em SIAM J. Control Optim.} {\bf 37}(1999), no. 6,
1568-1599.

\bibitem{Y1}P.F. Yao, Observability inequalities for the shallow
shell, {\em SIAM J. Control Optim.} {\bf 38}(2000), no. 6,
1729-1756.

\bibitem{Y2} P.F. Yao, Global smooth solutions for the quasilinear wave equation with
boundary dissipation, preprint, 2005.

\bibitem{YZ} J. Yong and X. Zhang, Exact controllability of the
heat equation with heperbolic memory kernel, {\em Control theory
of partial differentail equations}, 387-401, {\em Lect. Notes Pure
Appl. Math.}, {\bf 424} (2005).
\end{thebibliography}
\end{document}